\documentclass[11 pt,onecolumn]{IEEEtran}

\IEEEoverridecommandlockouts     

\usepackage{graphicx,algorithm,algorithmic,subfigure,latexsym,amsmath,amsfonts,amssymb,cite,mathrsfs}
\usepackage{epsfig,color}
\usepackage{epstopdf}

\newtheorem{lemma}{Lemma}[section]
\newtheorem{theorem}{Theorem}[section]
\newtheorem{corollary}{Corollary}[section]
\newtheorem{remark}{Remark}[section]
\newtheorem{definition}{Definition}[section]
\newtheorem{example}{Example}[section]

\newtheorem{assumption}{Assumption}[section]

\def\rrr#1\\{\par
\medskip\hbox{\vbox{\parindent=2em\hsize=6.12in
\hangindent=4em\hangafter=1#1}}}
\allowdisplaybreaks

\begin{document}



%

\title{\LARGE Semistability-Based Robust and Optimal Control Design for Network Systems}

\author{Qing Hui and Zhenyi Liu\\
\small{Control Science and Engineering Laboratory\\
Department of Mechanical Engineering \\
Texas Tech University \\
Lubbock, TX 79409-1021\\
{\tt\small (qing.hui@ttu.edu; zhenyi.liu@ttu.edu)}\\
Technical Report CSEL-07-14, July 2014}
\thanks{This work was supported by the Defense Threat Reduction Agency, Basic Research Award \#HDTRA1-10-1-0090 and Fundamental Research Award \#HDTRA1-13-1-0048, to Texas Tech University.}
}

\maketitle
\IEEEpeerreviewmaketitle
\pagestyle{empty}
\thispagestyle{empty}
\baselineskip 22pt

\begin{abstract}
In this report, we present a new Linear-Quadratic Semistabilizers (LQS) theory for
linear network systems. This new semistable $\mathcal{H}_{2}$ control framework is developed to address the robust and optimal semistable control issues of network systems while preserving network topology subject to white noise. Two new notions of semistabilizability and semicontrollability are introduced as a means to connecting semistability with the Lyapunov equation based technique. With these new notions, we first develop a semistable $\mathcal{H}_{2}$ control theory for network systems by exploiting the properties of semistability. A new series of necessary and sufficient conditions for semistability of the closed-loop system have been derived in terms of the Lyapunov equation. Based on these results, we propose a constrained optimization technique to solve the semistable $\mathcal{H}_{2}$ network-topology-preserving control design for network systems over an admissible set. Then optimization analysis and the development of numerical algorithms for the obtained constrained optimization problem are conducted. We establish the existence of optimal solutions for the obtained nonconvex optimization problem over some admissible set. Next, we propose a heuristic swarm optimization based numerical algorithm towards efficiently solving this nonconvex, nonlinear optimization problem. Finally, several numerical examples will be provided.
\end{abstract}

\section{Introduction}

Linear-Quadratic Semistabilizers (LQS) is a generalization of deterministic \textit{optimal semistable control} proposed in \cite{HH:IJC:2009,Hui:SCL:2011,HHC:JFI:2011,Hui:AUT:2011,Hui:JFI:2012} to the stochastic case in which
the stochasticity comes from two aspects: random distribution of initial conditions and stochastic $L_{2}$ noises due to sensor noise and exogenous disturbance. LQS also has a similar mathematical formulation as the stochastic LQR control while their differences are apparent. First, LQR control guarantees asymptotic stability of the closed-loop system, that is, state convergence to the origin, while LQS control assures \textit{semistability} of the closed-loop system. As discussed in \cite{HHB:TAC:2008}, semistability is the     
property that the limiting state is determined not only by system dynamics, but also by initial conditions, and hence, is not fixed \textit{a priori}. Thus, LQS is essentially linear-quadratic control with \textit{nondeterministic} steady-state regulation. Next, LQR assures one unique feedback controller while, on the contrary, LQS allows multiple feedback controllers. Actually LQS may have infinite many solutions due to the initial results for the deterministic case developed in \cite{HH:IJC:2009,Hui:SCL:2011,HHC:JFI:2011,Hui:AUT:2011,Hui:JFI:2012}. Thus, this gap between the nonuniqueness feature of LQS and the uniqueness property of LQR leads to a huge difference in theoretical analysis and control applications of LQS. In this report, we initiate the first, systemic control-theoretic framework for the LQS theory based on many preliminary results \cite{HH:IJC:2009,Hui:SCL:2011,HHC:JFI:2011,Hui:AUT:2011,Hui:JFI:2012} that we have developed before. We hope this report will lead to research attention on LQS and its related control problems for network systems. 

The motivation of this report can be illustrated by the following three examples from mechanical engineering, network sciences, and electrical systems, respectively. The first motivating example is a simple mass-damper system.
For mechanical systems, semistability characterizes the motion of a rigid body subject to damping and in the absence of stiffness \cite{BB:JVA:1995}. Such a damped rigid body converges to rest at a Lyapunov stable equilibrium position that is determined by the initial conditions. Specifically, we consider the mass-damper system given by the equation of motion $m\ddot{q}(t)+c\dot{q}(t)=0$ with the initial condition $(q(0),\dot{q}(0))=(q_{0},\dot{q}_{0})$, $t\geq0$, where $q(t)\in\mathbb{R}$ denotes the position, $m>0$ denotes the mass, and $c>0$ denotes the viscous damping constant. It is straightforward to verify that the system state $(q(t),\dot{q}(t))$ converges to $(q_{\infty},0)$, where $q_{\infty}=\lim_{t\to\infty}q(t)=q_{0}+(m/c)\dot{q}_{0}$ is the limiting position determined by both the system parameters $m,c$ and the initial condition. Our control problem here is how to choose the ratio $m/c$ such that when the initial condition $(q(0),\dot{q}(0))$ is randomly distributed, the variance of the weighted convergence error $\lim_{t\to\infty}\mathbb{E}[(q(t)-q_{\infty})^{2}+k(\dot{q}(t))^{2}]$ is minimized, where $\mathbb{E}$ denotes the expectation operator and $k>0$ is a weight constant. The physical meaning for this control design is that we try to attenuate the stochastic effect as much as we can so that the system's response is robust to such an effect. Note that $q_{\infty}$ is not a fixed point and cannot be predetermined. 

The second motivating example is the consensus problem with imperfect information. In many applications involving multiagent coordination and distributed computation \cite{HHB:ACC:2008,HHB:TAC:2008,Hui:MTNS:2010,Hui:TAC:2011}, groups of agents
are required to agree on certain quantities of interest. In
particular, it is important to develop information consensus
protocols for networks of dynamic agents wherein a unique feature of
the closed-loop dynamics under any control algorithm that achieves
consensus is the existence of a continuum of equilibria representing
a state of agreement or \textit{consensus}. Technically, we consider a multiagent network system consisting of $n$ agents, whose dynamics are $\dot x_i = u_i$, where $x_i \in
\mathbb{R}$ is the state and $u_i$ is the control input. The linear consensus protocol can be designed as
\begin{eqnarray}
u_i = \sum_{j \in \mathcal{N}_i} a_{ij}(x_j - x_i),
\end{eqnarray}
where $a_{ij}>0$ and $\mathcal{N}_i$ denotes the set of all the \emph{neighbors} for agent $i$. The equivalent matrix form of this consensus protocol is given by
$\dot{x}(t)=-Lx(t)$, $x(0)=x_{0}$, $t\geq0$, where $L$ denotes a weighted \textit{Laplacian} matrix. It can be shown that if the underlying graph for $L$ is chosen to be undirected and connected \cite{HHB:TAC:2008}, then $\lim_{t\to\infty}x(t)=x_{\infty}=c[1,\ldots,1]^{\rm{T}}$, where $c$ is determined by the initial condition $x(0)$. Now we consider the case where $x(t)$ is imperfect and corrupted by some noise. In this case, the previous consensus protocol becomes
\begin{eqnarray}\label{eq:xdyna} 
\dot{x}(t)=-L(x(t)+w(t)), 
\end{eqnarray} where $w(t)$ denotes the standard Gaussian white noise. Now the design problem becomes how to choose $L$ such that the variance of the weighted convergence error $\lim_{t\to\infty}\mathbb{E}[(x(t)-x_{\infty})^{\rm{T}}R(x(t)-x_{\infty})]$ is minimized, where $R$ is symmetric and positive semidefinite.  

The third motivating example comes from \cite{SDD:CDC:2007}. The physical meaning of LQG heat engines is to seek an ideal current source to extract maximal amount of energy from resistors with Johnson-Nyquist noise. From a control-theoretic point of view, we aim to design an LQG optimal controller to extract work from heat flows for a heated lossless system. Specifically, consider the linear system given by $\dot{x}(t)=Jx(t)+Bu(t)$ and $y(t)=B^{\rm{T}}x(t)$ where $J$ is skew-symmetric and $(J,B)$ is controllable. The control aim here is to design a controller $u(t)=Kx(t)+Dw(t)$ to maximize the expected work extracted from this linear system and to achieve ``\textit{thermal equilibrium}", where the expected work is given by $-\int_{0}^{t}(\mathbb{E}[(x(s)-\mathbb{E}x(s))^{\rm{T}}Bu(s)]+u^{\rm{T}}(s)Ru(s)){\rm{d}}s$. This thermal equilibrium implies the temperature equipartition which is a semistable status shown by \cite{Hui:ACC:2011}. 

In the first part of this report, we generalize the above ideas of designing optimal mass-damper systems, consensus protocols, and LQG heat engines with an $\mathcal{H}_{2}$ performance to develop a semistable $\mathcal{H}_{2}$ network-topology-preserving control design framework for network systems by using two new notions of semistabilizability/semidetectability originating from \cite{HL:Allerton:2011}, and semicontrollability/semiobservability originating from \cite{HHC:JFI:2011}. Our proposed LQS framework not only works for the consensus problems and electrical circuits, but also applies to a wide range of applications such as chemical reaction systems, biomedical systems, and flight control systems in which semistability is the appropriate notion of stability \cite{BB:ACC:1999,HHB:TAC:2008,HCH:2009}. Although $\mathcal{H}_{2}$ optimal semistable control for linear systems has been discussed in \cite{HHC:JFI:2011,Hui:SCL:2011}, only the deterministic, free-network-topology case has been considered. However, it is much more involved and important to develop a stochastic $\mathcal{H}_{2}$ semistable network-topology-preserving control framework for network systems to address the robustness issue of such a design. Hence, this report takes on this issue to try to fill in the theoretical gap of semistable $\mathcal{H}_{2}$ theory for network systems.

The aim of this new type of control is to address emerging stability issue arising from network systems together with its optimal performance in a noisy networked environment. In particular, we want to design a distributed state feedback controller for a class of stochastic linear networks so that it can guarantee that the closed-loop network system is semistable and its $\mathcal{H}_{2}$ performance cost functional is minimized. There are many important applications for which semistability is the most appropriate stability property of interest. A classical example is the synchronization of multiple weakly coupled oscillators to a common frequency. Recently, significant results have been obtained on semistability in consensus problems for networked agents \cite{CHHR:SCL:2008,HHB:TAC:2008,HH:AUT:2008,HH:IJC:2009,Hui:TAC:2011,Hui:TAC:2012}. An example of such a problem is for a group of networked autonomous vehicles to converge to a common heading, and for the network to respond to a small perturbation with only a corresponding small change to the common heading. Semistability arises naturally in the context of system thermodynamics. System thermodynamics uses basic thermodynamic principles to design controllers that guarantee a high degree of robustness as indicated in \cite{HCN:2005}. Other recent results in semistability theory can be found in \cite{BB:ACC:1999,BB:SICON:2003,BB:ACC:2003}.

In the second part of this report, we address theoretic optimization analysis for a nonconvex optimization problem which characterizes the proposed LQS design derived from the first part of the report. The existence of optimal solutions will be established for this optimization formulation over some admissible set by using a series of important results. To our best knowledge, this is the first attempt to answering this open question regarding the existence of optimal solutions since many optimization-based forms of semistable control have been developed in recent years \cite{HH:IJC:2009,Hui:SCL:2011,HHC:JFI:2011,Hui:AUT:2011,Hui:JFI:2012}. Furthermore, inspired by the resemblance between swarm behaviors such as particle swarm \cite{KE:ICNN:1995} and semistable control for state equipartitioning \cite{CHHR:SCL:2008}, next we will propose to use a heuristic swarm optimization based numerical algorithm to solve the proposed nonconvex optimization problem. Finally, several numerical examples will be provided. A preliminary version of this report has appeared in \cite{HL:CDC:2012b}.

\section{Problem Formulation for Network Systems}\label{PFNS}

To generalize the three examples including optimal consensus problems presented in the Introduction, we use the semistability theory developed for linear systems \cite{HHC:JFI:2011} to formulate a robust and optimal semistable control problem for general network systems.
Before we proceed, we need to use graph notations to represent
network topologies for general network systems. Specifically, let $\mathcal{G}=(\mathcal{V},\mathcal{E},\mathcal{A})$
be a \textit{directed graph} (or digraph) with the set of \textit{nodes}
(or vertices) $\mathcal{V}=\{1,\dots,n\}$, the set of \textit{edges}
$\mathcal{E}\subseteq\mathcal{V}\times\mathcal{V}$ involving a set of ordered pairs, and an \textit{adjacency} matrix
$\mathcal{A}\in\mathbb{R}^{n\times n}$ such that $\mathcal{A}_{(i,j)}=1$,
$i,j=1,\dots,n$, if $(j,i)\in\mathcal{E}$, and $0$ otherwise.  Moreover, we assume that $\mathcal{A}_{(i,i)}=1$ for all
$i\in\mathcal{V}$. A \emph{graph} or \textit{undirected graph}
$\mathcal{G}$ associated with the adjacency matrix
$\mathcal{A}\in\mathbb{R}^{n\times n}$ is a directed graph for which
the \emph{arc set} is symmetric, that is,
$\mathcal{A}=\mathcal{A}^{\rm{T}}$. A graph $\mathcal{G}$ is
\textit{balanced} if
$\sum_{j=1}^{n}\mathcal{A}_{(i,j)}=\sum_{j=1}^{n}\mathcal{A}_{(j,i)}$
for all $i=1,\ldots,n$. 

For a graph or a set of graphs $\mathcal{G}$, we consider a dynamical network in which node $i$ has a \emph{value}{\footnote{It is straightforward to extend $x_{i}(t)$, $u_{i}(t)$, and $w_{i}(t)$ from the scalar case $\mathbb{R}$ to the vector case $\mathbb{R}^{q}$.}} given by $x_{i}(t)\in\mathbb{R}$ and a \textit{control input} $u_{i}(t)\in\mathbb{R}$, $i=1,\ldots,n$, at time $t\in\mathbb{R}$. We assume that this dynamical network is affected by $n$ independent, additive Gaussian white noises $w_{i}(t)$, $i=1,\ldots,n$. The evolution process for this dynamical network is governed by the following interconnected differential equations
\begin{eqnarray}\label{NS}
\dot{x}_{i}(t)&=&\mathcal{A}_{(i,i)}[a_{ii}x_{i}(t)+d_{ii}w_{i}(t)]+\sum_{j=1,j\neq i}^{n}\mathcal{A}_{(i,j)}[a_{ji}x_{j}(t)+d_{ji}w_{j}(t)\nonumber\\
&&-a_{ij}x_{i}(t)-d_{ij}w_{i}(t)]+u_{i}(t),\quad x_{i}(0)=x_{i0},\quad t\geq0,\quad i=1,\ldots,n,
\end{eqnarray} where $a_{ij},d_{ij}\in\mathbb{R}$ denote the constant weights for node values and white noises associated with the graph $\mathcal{G}$, respectively. We assume that $a_{ij},d_{ij}\equiv0$ if $\mathcal{A}_{(i,j)}=0$. The compact form of (\ref{NS}) can be written as 
\begin{eqnarray}\label{plant_SLQR}
\dot{x}(t)=A(\mathcal{G})x(t)+u(t)+D(\mathcal{G})w(t),\quad x(0)=x_{0},\quad t\geq0,\label{state}
\end{eqnarray} where $x(t)=[x_{1}(t),\ldots,x_{n}(t)]^{\rm{T}}\in\mathbb{R}^{n}$ is the system state vector, $u(t)=[u_{1}(t),\ldots,u_{n}(t)]^{\rm{T}}\in\mathbb{R}^{n}$ is the control input vector,
$w(t)=[w_{1}(t),\ldots,w_{n}(t)]^{\rm{T}}\in\mathbb{R}^{n}$ is the standard Gaussian white noise vector, $A=A(\mathcal{G})\in\mathbb{R}^{n\times n}$, and $D=D(\mathcal{G})\in\mathbb{R}^{n\times n}$. 

Note that (\ref{NS}) represents a wide spectrum of network systems in many disciplines. For example, the \textit{mass balance equation} for compartmental systems in \cite{HCH:2009} is indeed the form of (\ref{NS}). The linear inhibitory-excitatory neuronal network in \cite{HHB:NAHS:2011} has the same form as (\ref{NS}). Of course, the consensus protocol (\ref{eq:xdyna}) can also be viewed as (\ref{NS}). Thus, the dynamical network given by the form (\ref{NS}) offers a good foundation to conduct control synthesis of linear network systems.

The control aim here is to design a network state feedback controller given by
\begin{eqnarray}\label{control_G}
u_{i}(t)=\mathcal{A}_{(i,i)}k_{ii}x_{i}(t)+\sum_{j=1,j\neq i}^{n}\mathcal{A}_{(i,j)}[k_{ji}x_{j}(t)-k_{ij}x_{i}(t)],
\end{eqnarray} where $k_{ij}\in\mathbb{R}$ satisfies that $k_{ij}\equiv0$ if $\mathcal{A}_{(i,j)}=0$, or, equivalently in vector form,
\begin{eqnarray}\label{control}
u(t)=K(\mathcal{G})x(t),
\end{eqnarray} where $K=K(\mathcal{G})\in\mathbb{R}^{n\times n}$, such that the following design criteria are satisfied:
\begin{itemize}
\item[$i$)] The closed-loop system (\ref{state}) and (\ref{control}) is semistable, i.e., $\tilde{A}=A+K$ is semistable. 
\item[$ii$)] The performance functional 
\begin{eqnarray}\label{cost_SLQR}
J(K)=\lim_{t\to\infty}\frac{1}{t}\mathbb{E}\Big\{\int_{0}^{t}\Big[(x(s)-x_{\infty})^{\rm{T}}R_{1}(x(s)-x_{\infty})+(u(s)-u_{\infty})^{\rm{T}}R_{2}(u(s)-u_{\infty})\Big]{\rm{d}}s\Big\}
\end{eqnarray} is minimized, where $R_{1}=E_{1}^{\rm{T}}E_{1}$, $R_{2}=E_{2}^{\rm{T}}E_{2}$, $E_{1}^{\rm{T}}E_{2}=0$, $E_{1},E_{2}\in\mathbb{R}^{q\times n}$, $x_{\infty}=\lim_{t\to\infty}\mathbb{E}[x(t)]$, and $u_{\infty}=Kx_{\infty}$.
\end{itemize} The state feedback controllers satisfying Conditions $i$) and $ii$) are called \textit{linear-quadratic semistabilizers}.  Note that (\ref{control_G}) has the same network topology as (\ref{NS}). Furthermore, the closed-loop dynamics matrix $\tilde{A}$ has the same network topology as (\ref{NS}). Hence, (\ref{control_G}) is a fixed-structure controller which preserves the network topology $\mathcal{G}$. 

Recall from \cite{Bernstein:2009} that a matrix $A\in\mathbb{R}^{n\times n}$ is \textit{semistable} if and only if $\lim_{t\to\infty}e^{At}$ exists. Furthermore, if $A$ is semistable, then the index of $A$ is zero or one, and thus $A$ is group invertible. The \textit{group inverse} $A^{\#}$ of $A$ is a special case of the Drazin inverse $A^{\rm{D}}$ in the case in which $A$ has index zero or one \cite{Bernstein:2009}. In this case, $\lim_{t\to\infty}e^{At}=I_{q}-AA^{\#}$ \cite{Bernstein:2009}. Clearly, this control problem can be viewed as a stochastic version of semistable LQR control. However, due to the possibility of singularity of $\tilde{A}$, the analysis of such a control problem is much more involved than the standard LQR theory. In fact, the three examples in \cite{Hui:SCL:2011} show that this semistable $\mathcal{H}_{2}$ control problem could have none, one, or infinitely many solutions under the standard assumptions from the LQR theory. Hence, there is a big difference between the standard $\mathcal{H}_{2}$ control theory and the proposed stochastic optimal semistable control problem here. 

One may argue that we can resort to the
existing theory of Lyapunov and Riccati equations
to tackle semistability and synthesis
of linear systems with additive noise. For instance,
checking semistability via Lyapunov equations:
it is well known that if there exists a solution to
the Lyapunov equation then the modes (in the setup of this report, the pair $(A^{\rm{T}},D^{\rm{T}})$) that are
observable are asymptotically stable. One simply has 
to check that the remaining modes (that is, 
non-observable modes) are semistable. Such a method relies on the Kalman decomposition of $(A^{\rm{T}},D^{\rm{T}})$ and this matrix transformation will destroy the graph topology $\mathcal{G}$ in $A$ and $D$. However, the control design needs to retain the graph topology in $K$ shown in (\ref{control_G}) since the topologically structured information is the only one available for feedback in many network systems. As an example, in sensor networks the relative sensing information is always used for control instead of absolute sensing information \cite{ZM:TAC:2011}. In this case, it is more realistic to use the form of (\ref{control_G}) to design controllers for sensor network systems. Thus, the matrix decomposition method which disrupts the network topology becomes very difficult to address this property required in control. 

\section{Mathematical Preliminaries}\label{MP}

Suppose $\tilde{A}$ is semistable. Then it follows from (\ref{cost_SLQR}) and Lemma~\ref{limit} in the Appendix that 
\begin{eqnarray}\label{cost_J}
J(K)&=&\lim_{t\to\infty}\mathbb{E}\left[\tilde{x}^{\rm{T}}(t)\tilde{R}\tilde{x}(t)\right],
\end{eqnarray} where $\tilde{x}(t)=x(t)-x_{\infty}$ and $\tilde{R}=R_{1}+K^{\rm{T}}R_{2}K$. Next, it follows from (\ref{plant_SLQR}) and (\ref{control}) that 
\begin{eqnarray}\label{Ext}
\frac{{\rm{d}}}{{\rm{d}}t}\mathbb{E}[x(t)]=\tilde{A}\mathbb{E}[x(t)],\quad t\geq0.
\end{eqnarray} Hence, $\mathbb{E}[x(t)]=e^{\tilde{A}t}\mathbb{E}[x(0)]$ and $x_{\infty}=\lim_{t\to\infty}\mathbb{E}[x(t)]=(I_{n}-\tilde{A}\tilde{A}^{\#})\mathbb{E}[x(0)]$. Note that $x_{\infty}\in\mathcal{N}(\tilde{A})$, where $\mathcal{N}(X)$ denotes the null space of $X$. Then it follows from (\ref{plant_SLQR}) and (\ref{control}) that 
\begin{eqnarray}
&&\dot{\tilde{x}}(t)=\tilde{A}\tilde{x}(t)+Dw(t),\label{xtilde}\\
&&\frac{{\rm{d}}}{{\rm{d}}t}\mathbb{E}[\tilde{x}(t)]=\tilde{A}\mathbb{E}[\tilde{x}(t)],\quad t\geq0,\label{1stM}
\end{eqnarray} where $\mathbb{E}[\tilde{x}(0)]=\tilde{A}\tilde{A}^{\#}\mathbb{E}[x(0)]$. Note that (\ref{Ext}) is the first moment equation for (\ref{plant_SLQR}) and (\ref{control}) while (\ref{1stM}) is the first moment equation for (\ref{xtilde}). To completely analyze (\ref{plant_SLQR}) and (\ref{control}) as well as (\ref{xtilde}), one needs to consider the second moment equations. Define the covariance matrices $Q(t)=\mathbb{E}[x(t)x^{\rm{T}}(t)]-\mathbb{E}[x(t)]\mathbb{E}[x^{\rm{T}}(t)]$ and $\tilde{Q}(t)=\mathbb{E}[\tilde{x}(t)\tilde{x}^{\rm{T}}(t)]-\mathbb{E}[\tilde{x}(t)]\mathbb{E}[\tilde{x}^{\rm{T}}(t)]$. We make the following standing assumption.

\begin{assumption}\label{A1}
$x(0)$ and $w(t)$ are independent for all $t\geq0$.
\end{assumption}

The following result asserts that $Q(t)$ and $\tilde{Q}(t)$ are the same. Hence, one can use either $Q(t)$ or $\tilde{Q}(t)$ to characterize the covariance matrix for the closed-loop system.

\begin{lemma}\label{lemma_Qt}
Assume that Assumption~\ref{A1} holds. Then $Q(t)$ satisfies the Lyapunov differential equation
\begin{eqnarray}\label{LE}
\dot{Q}(t)=\tilde{A}Q(t)+Q(t)\tilde{A}^{\rm{T}}+DD^{\rm{T}},\quad t\geq0,
\end{eqnarray} where $Q(0)=\mathbb{E}[x(0)x^{\rm{T}}(0)]-\mathbb{E}[x(0)]\mathbb{E}[x^{\rm{T}}(0)]$. Furthermore, (\ref{LE}) is equivalent to 
\begin{eqnarray}\label{Qt}
Q(t)=e^{\tilde{A}t}Q(0)e^{\tilde{A}^{\rm{T}}t}+\int_{0}^{t}e^{\tilde{A}s}DD^{\rm{T}}e^{\tilde{A}^{\rm{T}}s}{\rm{d}}s.
\end{eqnarray} Finally, $Q(t)=\tilde{Q}(t)\geq0$ for all $t\geq0$. 
\end{lemma}

\begin{IEEEproof}
Note that $\dot{x}(t)=\tilde{A}x(t)+Dw(t)$. Evaluating $\dot{Q}(t)$ yields
\begin{eqnarray}\label{dQt}
\dot{Q}(t)&=&\mathbb{E}[\dot{x}(t)x^{\mathrm{T}}(t)+x(t)\dot{x}^{\mathrm{T}}(t)]-\mathbb{E}[\dot{x}(t)]\mathbb{E}[x^{\mathrm{T}}(t)]-\mathbb{E}[x(t)]\mathbb{E}[\dot{x}^{\mathrm{T}}(t)]\nonumber\\
&=&\mathbb{E}[(\tilde{A}x(t)+Dw(t))x^{\mathrm{T}}(t)+x(t)(\tilde{A}x(t)+Dw(t))^{\mathrm{T}}]-\tilde{A}\mathbb{E}[x(t)]\mathbb{E}[x^{\mathrm{T}}(t)]-\mathbb{E}[x(t)]\mathbb{E}[x^{\mathrm{T}}(t)]\tilde{A}^{\mathrm{T}}\nonumber\\
&=&\mathbb{E}[\tilde{A}x(t)x^{\mathrm{T}}(t)]-\tilde{A}\mathbb{E}[x(t)]\mathbb{E}[x^{\mathrm{T}}(t)]+\mathbb{E}[x(t)x^{\mathrm{T}}(t)\tilde{A}^{\mathrm{T}}]-\mathbb{E}[x(t)]\mathbb{E}[x^{\mathrm{T}}(t)]\tilde{A}^{\mathrm{T}}\nonumber\\
&&+\mathbb{E}[Dw(t)x^{\mathrm{T}}(t)+x(t)w^{\mathrm{T}}(t)D^{\mathrm{T}}]\nonumber\\
&=&\tilde{A}Q(t)+Q(t)\tilde{A}^{\mathrm{T}}+\mathbb{E}\Big[Dw(t)\Big(e^{\tilde{A}t}x_{0}+\int_{0}^{t}e^{\tilde{A}(t-s)}Dw(s){\mathrm{d}}s\Big)^{\mathrm{T}}\nonumber\\
&&+\Big(e^{\tilde{A}t}x_{0}+\int_{0}^{t}e^{\tilde{A}(t-s)}Dw(s){\mathrm{d}}s\Big)w^{\mathrm{T}}(t)D^{\mathrm{T}}\Big]\nonumber\\
&=&\tilde{A}Q(t)+Q(t)\tilde{A}^{\mathrm{T}}+\mathbb{E}\Big[D\int_{0}^{t}w(t)w^{\mathrm{T}}(s)D^{\mathrm{T}}e^{\tilde{A}^{\mathrm{T}}(t-s)}{\mathrm{d}}s+\int_{0}^{t}e^{\tilde{A}(t-s)}Dw(s)w^{\mathrm{T}}(t){\mathrm{d}}sD^{\mathrm{T}}\Big]\nonumber\\
&=&\tilde{A}Q(t)+Q(t)\tilde{A}^{\mathrm{T}}+D\int_{0}^{t}\delta(t-s)D^{\mathrm{T}}e^{\tilde{A}^{\mathrm{T}}(t-s)}{\mathrm{d}}s+\int_{0}^{t}e^{\tilde{A}(t-s)}D\delta(s-t){\mathrm{d}}sD^{\mathrm{T}}\nonumber\\
&=&\tilde{A}Q(t)+Q(t)\tilde{A}^{\mathrm{T}}+\frac{1}{2}DD^{\mathrm{T}}+\frac{1}{2}DD^{\mathrm{T}}\nonumber\\
&=&\tilde{A}Q(t)+Q(t)\tilde{A}^{\mathrm{T}}+DD^{\mathrm{T}},
\end{eqnarray} where we used the fact that $\mathbb{E}[w(t)w^{\mathrm{T}}(s)]=\delta(t-s)$ and $\delta(t)$ denotes the Dirac function. Hence, (\ref{LE}) holds.

Next, to show (\ref{Qt}), rewriting (\ref{LE}) as
\begin{eqnarray*}
{\mathrm{vec}}\,\dot{Q}(t)&=&{\mathrm{vec}}\,\tilde{A}Q(t)+{\mathrm{vec}}\,Q(t)\tilde{A}^{\mathrm{T}}+{\mathrm{vec}}\,DD^{\mathrm{T}}\nonumber\\
&=&(\tilde{A}\oplus\tilde{A}){\mathrm{vec}}\,Q(t)+{\mathrm{vec}}\,DD^{\mathrm{T}}
\end{eqnarray*} Now using the Lagrange's formula yields
\begin{eqnarray*}
{\mathrm{vec}}\,Q(t)=e^{(\tilde{A}\oplus\tilde{A})t}{\mathrm{vec}}\,Q(0)+\int_{0}^{t}e^{(\tilde{A}\oplus\tilde{A})(t-s)}{\mathrm{vec}}\,DD^{\mathrm{T}}{\mathrm{d}}s
\end{eqnarray*} or, equivalently, by changing the variable of integration
\begin{eqnarray*}
{\mathrm{vec}}\,Q(t)&=&e^{(\tilde{A}\oplus\tilde{A})t}{\mathrm{vec}}\,Q(0)+\int_{0}^{t}e^{(\tilde{A}\oplus\tilde{A})s}{\mathrm{vec}}\,DD^{\mathrm{T}}{\mathrm{d}}s\nonumber\\
&=&e^{\tilde{A}t}\otimes e^{\tilde{A}t}{\mathrm{vec}}\,Q(0)+\int_{0}^{t}e^{\tilde{A}s}\otimes e^{\tilde{A}s}{\mathrm{vec}}\,DD^{\mathrm{T}}{\mathrm{d}}s\nonumber\\
&=&{\mathrm{vec}}\,e^{\tilde{A}t}Q(0)e^{\tilde{A}^{\mathrm{T}}t}+\int_{0}^{t}\mathrm{vec}\,e^{\tilde{A}s}DD^{\mathrm{T}}e^{\tilde{A}^{\mathrm{T}}s}{\mathrm{d}}s,
\end{eqnarray*} which implies (\ref{Qt}). Finally, using the similar arguments as in (\ref{dQt}), it follows that $\dot{\tilde{Q}}(t)=\tilde{A}\tilde{Q}(t)+\tilde{Q}(t)\tilde{A}^{\mathrm{T}}+DD^{\mathrm{T}}$. Thus, $\dot{\tilde{Q}}(t)=e^{\tilde{A}t}\tilde{Q}(0)e^{\tilde{A}^{\rm{T}}t}+\int_{0}^{t}e^{\tilde{A}s}DD^{\rm{T}}e^{\tilde{A}^{\rm{T}}s}{\rm{d}}s$ for all $t\geq0$. Now since $\tilde{Q}(0)=\mathbb{E}[\tilde{x}(0)\tilde{x}^{\mathrm{T}}(0)]-\mathbb{E}[\tilde{x}(0)]\mathbb{E}[\tilde{x}^{\mathrm{T}}(0)]=\mathbb{E}[(x(0)-x_{\infty})(x(0)-x_{\infty})^{\mathrm{T}}]-\mathbb{E}[x(0)-x_{\infty}]\mathbb{E}[(x(0)-x_{\infty})^{\mathrm{T}}]=\mathbb{E}[x(0)x^{\mathrm{T}}(0)]-\mathbb{E}[x(0)]\mathbb{E}[\\x^{\mathrm{T}}(0)]=Q(0)$, it follows that $Q(t)=\tilde{Q}(t)$ for all $t\geq0$.
\end{IEEEproof}

Now we have the following necessary and sufficient condition on the convergence of an integral in (\ref{Qt}).

\begin{lemma}\label{lemma_int}
Assume that Assumption~\ref{A1} holds. If $\tilde{A}$ is semistable, then $\hat{Q}\triangleq\lim_{t\to\infty}\int_{0}^{t}e^{\tilde{A}s}DD^{\rm{T}}e^{\tilde{A}^{\rm{T}}s}{\rm{d}}s$ exists if and only if $\mathcal{N}(\tilde{A}^{\rm{T}})\subseteq\mathcal{N}(D^{\rm{T}})$, where $\mathcal{N}(X)$ denotes the null space of $X$. Furthermore, if $\mathcal{N}(\tilde{A}^{\rm{T}})\subseteq\mathcal{N}(D^{\rm{T}})$, then $Q\triangleq\lim_{t\to\infty}Q(t)$ is given by 
\begin{eqnarray}\label{Q}
Q&=&\lim_{t\to\infty}\mathbb{E}[\tilde{x}(t)\tilde{x}^{\rm{T}}(t)]\nonumber\\
&=&(I_{n}-\tilde{A}\tilde{A}^{\#})Q(0)(I_{n}-\tilde{A}\tilde{A}^{\#})^{\rm{T}}+\hat{Q}.
\end{eqnarray} 
\end{lemma}

\begin{IEEEproof}
Note that since $\tilde{A}$ is
semistable, it follows from ${\mathrm{mspec}}\,(\tilde{A})={\mathrm{mspec}}\,(\tilde{A}^{\mathrm{T}})$ (see Proposition 4.4.5 of \cite[p.~263]{Bernstein:2009}) that $\tilde{A}^{\mathrm{T}}$ is semistable as well, where ${\mathrm{mspec}}\,(\tilde{A})$ denotes the multispectrum of $\tilde{A}$. Hence, it follows that either $\tilde{A}^{\mathrm{T}}$ is Hurwitz or there exists an
invertible matrix $S\in\mathbb{R}^{n\times n}$ such that $\tilde{A}^{\mathrm{T}}=S\left[
                                                               \begin{array}{cc}
                                                                 J & 0 \\
                                                                 0 & 0 \\
                                                               \end{array}
                                                             \right]S^{-1}
$, where $J\in\mathbb{R}^{r\times r}$, $r={\rm{rank}}\,\tilde{A}^{\mathrm{T}}$, and $J$
is Hurwitz. Now, if $\tilde{A}^{\mathrm{T}}$ is Hurwitz, then $\mathcal{N}(\tilde{A}^{\mathrm{T}})=\{0\}\subseteq\mathcal{N}(D^{\mathrm{T}})$.

Alternatively, if $\tilde{A}^{\mathrm{T}}$ is not Hurwitz, then
$\mathcal{N}(\tilde{A}^{\mathrm{T}})=\left\{x\in\mathbb{R}^{n}:x=S[0_{1\times
r},y^{\rm{T}}]^{\rm{T}},y\in\mathbb{R}^{n-r}\right\}$.
Now,
\begin{eqnarray}\label{int_1}
\int_{0}^{\infty}e^{\tilde{A}t}DD^{\mathrm{T}}e^{\tilde{A}^{\mathrm{T}}t}{\rm{d}}t&=&S^{-{\rm{T}}}\int_{0}^{\infty}e^{\hat{J}t}DD^{\mathrm{T}}e^{\hat{J}t}{\rm{d}}tS\nonumber\\
&=&S^{-{\rm{T}}}\int_{0}^{\infty}\left[
                                       \begin{array}{cc}
                                         e^{J^{\rm{T}}t}\hat{R}_{1}e^{Jt} & e^{J^{\rm{T}}t}\hat{R}_{12} \\
                                         \hat{R}_{12}^{\rm{T}}e^{Jt} & \hat{R}_{2} \\
                                       \end{array}
                                     \right]
{\rm{d}}tS,
\end{eqnarray} where
$\hat{J}=\left[
          \begin{array}{cc}
            J & 0 \\
            0 & 0 \\
          \end{array}
        \right]$ and $\hat{R}=S^{\rm{T}}DD^{\mathrm{T}}S=\left[
                                             \begin{array}{cc}
                                               \hat{R}_{1} & \hat{R}_{12} \\
                                               \hat{R}_{12}^{\rm{T}} & \hat{R}_{2} \\
                                             \end{array}
                                           \right]$.
Next, it follows from (\ref{int_1}) that
$\int_{0}^{\infty}e^{\tilde{A}t}DD^{\mathrm{T}}e^{\tilde{A}^{\mathrm{T}}t}{\rm{d}}t$ exists 
if and only if $\hat{R}_{2}=0$ or, equivalently,
$[0_{1\times r},y^{\rm{T}}]\hat{R}[0_{1\times
r},y^{\rm{T}}]^{\rm{T}}=0$, $y\in\mathbb{R}^{n-r}$,
which is further equivalent to $x^{\rm{T}}DD^{\mathrm{T}}x=0$,
$x\in\mathcal{N}(\tilde{A}^{\mathrm{T}})$. Hence, $\mathcal{N}(\tilde{A}^{\mathrm{T}})\subseteq\mathcal{N}(D^{\mathrm{T}})$.

Finally, the proof of $\mathcal{N}(\tilde{A}^{\mathrm{T}})\subseteq\mathcal{N}(D^{\mathrm{T}})$ implies the existence of $\int_{0}^{\infty}e^{\tilde{A}t}DD^{\mathrm{T}}e^{\tilde{A}^{\mathrm{T}}t}{\rm{d}}t$ is immediate by reversing
the steps of the proof given above.

The second part of the result is a direct consequence of Lemma~\ref{lemma_Qt} and $\lim_{t\to\infty}e^{\tilde{A}t}=I_{n}-\tilde{A}\tilde{A}^{\#}$.
\end{IEEEproof}

Thus, under Assumption~\ref{A1}, if $\tilde{A}$ is semistable and $\mathcal{N}(\tilde{A}^{\rm{T}})\subseteq\mathcal{N}(D^{\rm{T}})$, then (\ref{cost_J}) can be rewritten as
\begin{eqnarray}\label{JKQ}
J(K)&=&{\rm{tr}}\,Q\tilde{R}\nonumber\\
&=&{\rm{tr}}\,(I_{n}-\tilde{A}\tilde{A}^{\#})Q(0)(I_{n}-\tilde{A}\tilde{A}^{\#})^{\rm{T}}\tilde{R}+{\rm{tr}}\,\hat{Q}\tilde{R}.
\end{eqnarray}
Clearly if $\mathbb{E}[x(0)x^{\rm{T}}(0)]=\mathbb{E}[x(0)]\mathbb{E}[x^{\rm{T}}(0)]$, then $Q(0)=0$, and hence, $J(K)={\rm{tr}}\,\hat{Q}\tilde{R}$, where ${\rm{tr}}\,X$ denotes the trace of $X$. One of the sufficient conditions to guarantee this scenario is that $x(0)$ is deterministic. However, here we consider the general case where $x(0)$ is not necessarily deterministic. Without loss of generality, we make the following assumption on $x(0)$.

\begin{assumption}\label{A2}
$x(0)$ is a random variable having a covariance matrix $V$, that is, $\mathbb{E}[x(0)x^{\rm{T}}(0)]-\mathbb{E}[x(0)x^{\rm{T}}(0)]\\=V$. 
\end{assumption}

Note that it follows from (\ref{Q}) that $Q$ has two parts: $(I_{n}-\tilde{A}\tilde{A}^{\#})Q(0)(I_{n}-\tilde{A}\tilde{A}^{\#})^{\rm{T}}$ and $\hat{Q}$. Hence, to minimize $J(K)$, one has to minimize two cost functionals associated with both terms simultaneously. This is not a good strategy from the optimization point of view. To combine two separate forms in (\ref{Q}) into a compact form, the following result is the key.

\begin{lemma}\label{lemma_JWV}
Assume that Assumptions \ref{A1} and \ref{A2} hold. Furthermore, assume $\mathcal{N}(\tilde{A}^{\rm{T}})\subseteq\mathcal{N}(D^{\rm{T}})$. If $\tilde{A}$ is semistable, then
\begin{eqnarray}\label{JWV}
J(K)&=&{\rm{tr}}\,(W+V)\tilde{R},
\end{eqnarray} where
\begin{eqnarray}\label{Wmatrix}
W=\int_{0}^{\infty}e^{\tilde{A}s}[\tilde{A}V+V\tilde{A}^{\rm{T}}+DD^{\rm{T}}]e^{\tilde{A}^{\rm{T}}s}{\rm{d}}s.
\end{eqnarray}
\end{lemma}

\begin{IEEEproof}
Using the fact that $\int_{0}^{\infty}\frac{{\rm{d}}}{{\rm{d}}s}(e^{\tilde{A}s}Ve^{\tilde{A}^{{\rm{T}}}s}){\rm{d}}s=e^{\tilde{A}s}Ve^{\tilde{A}^{{\rm{T}}}s}|_{s=0}^{s=\infty}$, we have
\begin{eqnarray}\label{Vmatrix}
\int_{0}^{\infty}e^{\tilde{A}s}[\tilde{A}V+V\tilde{A}^{\rm{T}}]e^{\tilde{A}^{\rm{T}}s}{\rm{d}}s=(I_{n}-\tilde{A}\tilde{A}^{\#})V(I_{n}-\tilde{A}\tilde{A}^{\#})^{\rm{T}}-V.
\end{eqnarray} Now the result is immediate. 
\end{IEEEproof}

It is important to note that unlike the standard LQR theory, here $Q$ given by (\ref{Q}) does not satisfy the standard Lyapunov equation $0=\tilde{A}Q+Q\tilde{A}^{\rm{T}}+DD^{\rm{T}}$ due to the fact that $\lim_{t\to\infty}\dot{Q}(t)=(I_{n}-\tilde{A}\tilde{A}^{\#})DD^{\rm{T}}(I_{n}-\tilde{A}\tilde{A}^{\#})^{\rm{T}}$. Furthermore, $Q$ given by (\ref{Q}) is just positive semidefinite, not positive definite. If we define $\mathcal{S}=\{K\in\mathbb{R}^{m\times n}:A+BK\,\,{\rm{is}}\,\,{\rm{semistable}}\}$ and $\mathcal{S}^{+}=\{K\in\mathcal{S}:Q=Q^{\rm{T}}\geq0\}$, then $\mathcal{S}^{+}$ is a closed set, not an open set. Hence, the standard Lagrange multiplier techniques cannot be applied to this problem. 

\section{Semidetectability and Semiobservability}\label{SSH2}

\subsection{Semidetectability}\label{SSIAO}

In this subsection, we introduce two new notions of \textit{semistabilizability} and \textit{semidetectability} \cite{HL:Allerton:2011} to address semistability of $\tilde{A}$ by using the Lyapunov equation while keeping the topological structure $\mathcal{G}$. 

\begin{definition}\label{Defo}
Let $A\in\mathbb{R}^{n\times n}$, $B\in\mathbb{R}^{n\times l}$, and $C\in\mathbb{R}^{l\times n}$.
The pair $(A,B)$ is \textit{semistabilizable} if 
\begin{eqnarray}\label{rank_SS}
{\rm{rank}}\left[\begin{array}{cc}
B & \jmath\omega I_{n}-A\\
\end{array}\right]=n
\end{eqnarray} for every nonzero $\omega\in\mathbb{R}$, where ${\rm{rank}}$ denotes the rank of a matrix. The pair $(A,C)$ is \textit{semidetectable} if 
\begin{eqnarray}
{\rm{rank}}\left[\begin{array}{c}
C\\
\jmath\omega I_{n}-A\\
\end{array}\right]=n\label{rank}
\end{eqnarray} for every nonzero $\omega\in\mathbb{R}$.
\end{definition}

It is clear from Definition~\ref{Defo} that $(A,C)$ is semidetectable if and only if $(A^{\rm{T}},C^{\rm{T}})$ is semistabilizable. Furthermore, it is important to note that semistabilizability and semidetectability are \textit{different} from the standard notions of stabilizability and detectability in linear control theory.  Recall that $(A,B)$ is stabilizable if and only if ${\rm{rank}}\left[\begin{array}{cc}
B & \lambda I_{n}-A\\
\end{array}\right]=n$ for every $\lambda\in\mathcal{C}$ in the closed right half plane, and $(A,C)$ is detectable if and only if ${\rm{rank}}\left[\begin{array}{c}
C\\
\lambda I_{n}-A\\
\end{array}\right]=n$ for every $\lambda\in\mathbb{C}$ in the closed right half plane. Obviously if $(A,C)$ is detectable, then it is semidetectable, but not vice versa. Similar remarks hold for the notions of controllability and observability, that is, if $(A,C)$ is observable, then it is semidetectable, but not vice versa. Hence, semidetectability is a much weaker notion than observability and detectability. Since (\ref{rank}) is only concerned with the detectability of $(A,C)$ on the imaginary axis, we label this notion as semidetectability. 

\begin{remark}\label{rmk_RN}
It follows from Facts 2.11.1-2.11.3 of \cite[p.130-131]{Bernstein:2009} that (\ref{rank_SS}) and (\ref{rank}) are equivalent to 
\begin{eqnarray}\label{dim}
\dim[\mathcal{R}(\jmath\omega I_{n}-A)+\mathcal{R}(B)]=n
\end{eqnarray} and 
\begin{eqnarray}
\mathcal{N}(\jmath\omega I_{n}-A)\cap\mathcal{N}(C)=\{0\}
\end{eqnarray} respectively, where $\dim$ denotes the dimension and $\mathcal{R}(X)$ denotes the range space of $X$. As we can see later on, these two equivalent forms are aligned with (\ref{Range}) and (\ref{Null}), respectively. \hfill$\blacklozenge$
\end{remark}

\begin{example}
Consider 
$A=\left[\begin{array}{cc}
0 & 0 \\
0 & 0 \\
\end{array}\right]$ and $B=\left[\begin{array}{c}
0\\
1\\
\end{array}\right]$. Clearly $(A,B)$ is not stabilizable. However, one can check (\ref{rank_SS}) that $(A,B)$ is indeed  semistabilizable. Thus, for this pair of $(A,B)$, the standard $\mathcal{H}_{2}$ control problem is not well defined, but the semistable $\mathcal{H}_{2}$ control problem is well defined.  \hfill$\blacktriangle$
\end{example}

Similar to stabilizability, state feedback control does not change semistabilizability. The proof is identical to the case of stabilizability by use of the Sylvester's inequality for rank. 

\begin{lemma}
Let $A\in\mathbb{R}^{n\times n}$ and $B\in\mathbb{R}^{n\times l}$. If $(A,B)$ is semistabilizable, then $(A+BK,B)$ is semistabilizable.
\end{lemma}

\begin{IEEEproof}
Since $(A,B)$ is semistabilizable, it follows that ${\mathrm{rank}}\left[\begin{array}{cc}
B & \jmath\omega I_{n}-A \\
\end{array}\right]=n$ for all nonzero $\omega\in\mathbb{R}$. Hence, by Sylvester's inequality, for all nonzero $\omega\in\mathbb{R}$,
\begin{eqnarray}\label{SVL}
n&=&n+(l+n)-(l+n)\nonumber\\
&=&{\mathrm{rank}}\left[\begin{array}{cc}
B & \jmath\omega I_{n}-A \\
\end{array}\right]+{\mathrm{rank}}\left[\begin{array}{cc}
I_{l} & -K\\
0 & I_{n}\\
\end{array}\right]-(l+n)\nonumber\\
&\leq&{\mathrm{rank}}\left[\begin{array}{cc}
B & \jmath\omega I_{n}-A \\
\end{array}\right]\left[\begin{array}{cc}
I_{l} & -K\\
0 & I_{n}\\
\end{array}\right]\nonumber\\
&\leq&{\mathrm{rank}}\left[\begin{array}{cc}
B & \jmath\omega I_{n}-A \\
\end{array}\right]\nonumber\\
&=&n.
\end{eqnarray} Since $\left[\begin{array}{cc}
B & \jmath\omega I_{n}-A-BK \\
\end{array}\right]=\left[\begin{array}{cc}
B & \jmath\omega I_{n}-A \\
\end{array}\right]\left[\begin{array}{cc}
I_{l} & -K\\
0 & I_{n}\\
\end{array}\right]$, it follows from (\ref{SVL}) that 
\begin{eqnarray*}
{\mathrm{rank}}\left[\begin{array}{cc}
B & \jmath\omega I_{n}-A-BK \\
\end{array}\right]=n
\end{eqnarray*} for all noznero $\omega\in\mathbb{R}$. Thus, $(A+BK,B)$ is semistabilizable. 
\end{IEEEproof}

Next, we introduce the core technical result of this report by connecting semistability with the Lyapunov equation via semidetectability. Before we state it, we need the following lemma.

\begin{lemma}\label{lemma_BP}
Let $A\in\mathbb{R}^{n\times n}$. Then $A$ is semistable if and only if $\mathcal{N}(A)\cap\mathcal{R}(A)=\{0\}$ and ${\rm{spec}}\,(A)\subseteq\{s\in\mathbb{C}:s+s^{*}<0\}\cup\{0\}$, where $s^{*}$ denotes the complex conjugate and ${\rm{spec}}\,(X)$ denotes the spectrum of $X$.
\end{lemma}

\begin{IEEEproof}
If $A$ is semistable, then it follows from Definition 11.8.1 of \cite[p.~727]{Bernstein:2009} that ${\rm{spec}}\,(A)\subseteq\{s\in\mathbb{C}:s+s^{*}<0\}\cup\{0\}$ and either $A$ is Hurwitz or there exists an
invertible matrix $S\in\mathbb{R}^{n\times n}$ such that $A=S\left[
                                                               \begin{array}{cc}
                                                                 J & 0 \\
                                                                 0 & 0 \\
                                                               \end{array}
                                                             \right]S^{-1}
$, where $J\in\mathbb{R}^{r\times r}$, $r={\rm{rank}}\,A$, and $J$
is Hurwitz. Now, if $A$ is Hurwitz, then $\mathcal{N}(A)=\{0\}=\mathcal{N}(A)\cap\mathcal{R}(A)$.

Alternatively, if $A$ is not Hurwitz, then $\mathcal{N}(A)=\{S[0_{1\times r},y_{2}^{\mathrm{T}}]^{\mathrm{T}}:y_{2}\in\mathbb{R}^{n-r}\}$. In this case, for every $S[0_{1\times r},x_{2}^{\mathrm{T}}]^{\mathrm{T}}\in\mathcal{N}(A)\cap\mathcal{R}(A)$, there exists $z\in\mathbb{R}^{n}$ such that $S[0_{1\times r},x_{2}^{\mathrm{T}}]^{\mathrm{T}}=Az$. Hence, $S[0_{1\times r},x_{2}^{\mathrm{T}}]^{\mathrm{T}}=S\left[
                                                               \begin{array}{cc}
                                                                 J & 0 \\
                                                                 0 & 0 \\
                                                               \end{array}
                                                             \right]S^{-1}z$, that is,
                                                             \begin{eqnarray*}
                                                             \left[\begin{array}{c}
                                                             0 \\
                                                             x_{2}\\
                                                             \end{array}\right]=\left[
                                                                                                                            \begin{array}{cc}
                                                                                                                              J & 0 \\
                                                                                                                              0 & 0 \\
                                                                                                                            \end{array}
                                                                                                                          \right]S^{-1}z,
                                                             \end{eqnarray*} which implies that $x_{2}=0$. Thus, $\mathcal{N}(A)\cap\mathcal{R}(A)=\{0\}$.
                                                             
Conversely, we assume that $\mathcal{N}(A)\cap\mathcal{R}(A)=\{0\}$ and ${\rm{spec}}\,(A)\subseteq\{s\in\mathbb{C}:s+s^{*}<0\}\cup\{0\}$. If $A$ is nonsingular, then $A$ is Hurwitz, and hence, $A$ is semistable. Now we consider the case where $A$ is singular. In this case, let $x\in\mathcal{N}(A^{2})$. Then it follows from $A^{2}x=AAx=0$ that $Ax\in\mathcal{N}(A)$. Note that $Ax\in\mathcal{R}(A)$. It follows from $\mathcal{N}(A)\cap\mathcal{R}(A)=\{0\}$ that $Ax=0$, that is, $x\in\mathcal{N}(A)$. Hence, $\mathcal{N}(A^{2})\subseteq\mathcal{N}(A)$. On the other hand, clearly $\mathcal{N}(A)\subseteq\mathcal{N}(A^{2})$. Thus, $\mathcal{N}(A)=\mathcal{N}(A^{2})$. Then it follows from Proposition 5.5.8 of \cite[p.~323]{Bernstein:2009} that the eigenvalue $0$ of $A$ is semisimple. Finally, by Definition 11.8.1 of \cite[p.~727]{Bernstein:2009}, $A$ is semistable. 
\end{IEEEproof}

\begin{lemma}\label{CLx}
Let $A\in\mathbb{R}^{n\times n}$ and $C\in\mathbb{R}^{l\times n}$. If $A$ is semistable and $\mathcal{N}(A)\subseteq\mathcal{N}(C)$, then $CL=0$, where $L$ is given by 
\begin{eqnarray}\label{L_def}
L\triangleq I_{n}-AA^{\#}. 
\end{eqnarray}
\end{lemma}

\begin{IEEEproof}
First it follows from semistability of $A$ that $L$ is well defined. We show that $CLx\equiv0$ for any $x\in\mathbb{R}^{n}$. Suppose that there exists $x\neq0$ such that $CLx\neq0$. Then $Lx\not\in\mathcal{N}(C)$. Since $\mathcal{N}(A)\subseteq\mathcal{N}(C)$, it follows that $Lx\not\in\mathcal{N}(A)$. However, $ALx=A(I_{n}-AA^{\#})x=(A-AAA^{\#})x=0$, which implies that $Lx\in\mathcal{N}(A)$. This is a contradiction. Hence, $CLx\equiv0$ for any  $x\in\mathbb{R}^{n}$.    
\end{IEEEproof}

\begin{theorem}\label{thm_rank}
Let $A\in\mathbb{R}^{n\times n}$. Then $A$ is semistable if and only if there exist a positive integer $m$, an $m\times n$ matrix $C$, and an $n\times n$ matrix $P=P^{\rm{T}}>0$ such that the pair $(A,C)$ is semidetectable and 
\begin{eqnarray}\label{newLE}
0=A^{\rm{T}}P+PA+C^{\rm{T}}C.
\end{eqnarray}
\end{theorem}

\begin{IEEEproof}
If $A$ is semistable, then by Lemma~\ref{lemma_BP}, $\jmath\omega$, $\omega\neq0$, cannot be an eigenvalue of $A$ and hence, ${\rm{rank}}\,(\jmath\omega I_{n}-A)=n$ for every nonzero $\omega\in\mathbb{R}$. Thus, (\ref{rank}) holds for any $C\in\mathbb{R}^{m\times n}$ and any positive integer $m$. To prove the existence of a positive definite solution to (\ref{newLE}), we restrict $C$ to be the one satisfying $\mathcal{N}(A)\subseteq\mathcal{N}(C)$. For such a pair
$(A,C)$, let
\begin{eqnarray}\label{hatP_def}
\hat{P}=\int_{0}^{\infty}e^{A^{\rm{T}}t}C^{\rm{T}}Ce^{At}{\rm{d}}t.
\end{eqnarray} 
Then it follows from Proposition 2.2 of \cite{HHC:JFI:2011} that such a $\hat{P}$ is well defined. Clearly $\hat{P}=\hat{P}^{\rm{T}}\geq0$. Since $A$ is semistable, by
Lemma~\ref{lemma_BP}, $\mathcal{N}(A)\cap\mathcal{R}(A)=\{0\}$, it
follows from \cite[p.~119]{BP:79} that $A$ is group
invertible. Hence, it follows from (\ref{hatP_def}) that 
\begin{eqnarray}\label{APA}
A^{\rm{T}}\hat{P}+\hat{P}A&=&\int_{0}^{\infty}\frac{{\rm{d}}}{{\rm{d}}t}\left(e^{A^{\rm{T}}t}C^{\rm{T}}Ce^{At}\right){\rm{d}}t\nonumber\\
&=&(I_{n}-AA^{\#})^{\rm{T}}C^{\rm{T}}C(I_{n}-AA^{\#})-C^{\rm{T}}C.
\end{eqnarray}

Next, consider $L$ given by (\ref{L_def}) and note that
$L^{2}=L$. Hence, $L$ is the unique $n\times n$ matrix satisfying
$\mathcal{N}(L)=\mathcal{R}(A)$, $\mathcal{R}(L)=\mathcal{N}(A)$,
and $Lx=x$ for all $x\in\mathcal{N}(A)$. Now, define
\begin{eqnarray}\label{P_def}
P\triangleq\hat{P}+L^{\rm{T}}L.
\end{eqnarray} Clearly $P$ is symmetric. Next, we show that $P$ is positive definite. Consider
the function $V(x)=x^{\rm{T}}Px$, $x\in\mathbb{R}^{n}$. If
$V(x)=0$ for some $x\in\mathbb{R}^{n}$, then $\hat{P}x=0$ and $Lx=0$. It
follows from $i$) of Lemma 2.4 in \cite{HHC:JFI:2011} that
$x\in\mathcal{N}(A)$, and $Lx=0$ implies that $x\in\mathcal{R}(A)$.
Now, it follows from $ii$) of Lemma 2.4 in \cite{HHC:JFI:2011} that $x=0$.
Hence, $P$ is positive definite. Now since
$LA=A-AA^{\#}A=0$, it follows from (\ref{APA}), (\ref{P_def}), and Lemma~\ref{CLx} that
\begin{eqnarray}
A^{\rm{T}}P+PA+C^{\rm{T}}C&=&A^{\rm{T}}\hat{P}+\hat{P}A+C^{\rm{T}}C+A^{\rm{T}}L^{\rm{T}}L+L^{\rm{T}}LA\nonumber\\
&=&L^{\rm{T}}C^{\rm{T}}CL+(LA)^{\rm{T}}L+L^{\rm{T}}LA\nonumber\\
&=&0,
\end{eqnarray} which is (\ref{newLE}). Hence, the conclusion holds for $C$ satisfying $\mathcal{N}(A)\subseteq\mathcal{N}(C)$ (one obvious choice for $C$ is $A$).

The converse part follows from Corollary 11.9.1 in \cite[p.~730]{Bernstein:2009}. For completeness, we include it here. Note that $V(x)$ we defined above is positive definite. Moreover, $\dot{V}(x)=-x^{\rm{T}}C^{\rm{T}}Cx\leq0$, $x\in\mathbb{R}^{n}$, which implies that $A$ is Lyapunov stable. 

To show that $A$ is semistable, let $\jmath\omega\in{\rm{spec}}(A)$, where $\omega\in\mathbb{R}$ is nonzero, and let $x\in\mathbb{C}^{n}$ be an associated eigenevector. Then it follows from (\ref{newLE}) that 
\begin{eqnarray}
-x^{*}C^{\rm{T}}Cx&=&x^{*}(A^{\rm{T}}P+PA)x\nonumber\\
&=&x^{*}[(\jmath\omega I_{n}-A)^{*}P+P(\jmath\omega I_{n}-A)]x\nonumber\\
&=&0.
\end{eqnarray} Hence, $Cx=0$, and thus, 
${\rm{rank}}\left[\begin{array}{c}
C\\
\jmath\omega I_{n}-A\\
\end{array}\right]x=0$, which implies that $x=0$ by (\ref{rank}). This contradicts $x\neq0$. Consequently, $\jmath\omega\not\in{\rm{spec}}(A)$ for all nonzero $\omega\in\mathbb{R}$. Hence, the eigenvalue of $A$ is either a real/complex number with negative real part or is zero, and if the eigenvalue is zero, then it is semisimple. Therefore, $A$ is semistable.    
\end{IEEEproof}

\begin{remark}\label{rmk_sd}
It is important to note that if $A$ is semistable, then for \textit{every} semidetectable pair $(A,C)$, one cannot always find $P=P^{\rm{T}}>0$ such that (\ref{newLE}) holds. To see this, consider $A=\left[\begin{array}{cc}
-1 & 0 \\
0 & 0 \\
\end{array}\right]$ and $C=\left[\begin{array}{cc}
0 & 1 \\
\end{array}\right]$. Clearly $A$ is semistable. Moreover,  $\left[\begin{array}{c}
C\\
\jmath\omega I_{n}-A\\
\end{array}\right]=\left[\begin{array}{cc}
0 & 1 \\
\jmath\omega+1 & 0 \\
0 & \jmath\omega
\end{array}\right]$ is a full rank matrix for every $\omega\in\mathbb{R}$. Hence, $(A,C)$ is semidetectable. However, one can verify that there does not exist $P=P^{\rm{T}}>0$ such that (\ref{newLE}) holds. \hfill$\blacklozenge$
\end{remark}

Theorem~\ref{thm_rank} gives a necessary and sufficient condition for semistability by use of a Lyapunov equation. Since we do not assume the detectability of $(A,C)$, asymptotic stability cannot be reached. Instead, by assuming a weaker notion--semidetectability, we can guarantee semistability. This result inspires us to design semistable $\mathcal{H}_{2}$ and  $\mathcal{H}_{\infty}$ controllers using this new notion and optimization methods such as linear matrix inequalities. 

\begin{example}
Consider 
$A=\left[\begin{array}{cc}
-1 & 1 \\
1 & -1 \\
\end{array}\right]$ and $C=\left[\begin{array}{cc}
1 & -1 \\
\end{array}\right]$. Clearly, for $P$ given by
$P=\frac{1}{2}\left[\begin{array}{cc}
1 & 0 \\
0 & 1 \\
\end{array}\right]$, $P>0$ and (\ref{newLE}) holds. In this case, it can be verified that (\ref{rank}) holds for $(A,C)$, and hence, $(A,C)$ is semidetectable. Hence, it follows from Theorem~\ref{thm_rank} that $A$ is semistable. In fact, 
$\lim_{t\to\infty}e^{At}=\frac{1}{2}\left[\begin{array}{cc}
1 & 1 \\
1 & 1 \\
\end{array}\right]$. \hfill$\blacktriangle$
\end{example}

\subsection{Semiobservability}

Theorem~\ref{thm_rank} gives a necessary and sufficient condition for guaranteeing semistability of linear time-invariant systems by using a rank condition and the Lyapunov equation. However, as we mentioned in Remark~\ref{rmk_sd}, the statement of Theorem~\ref{thm_rank} cannot be generalized to the case where the pair $(A,C)$ is arbitrary. This poses a gap between semidetectability and detectability since the classic Lyapunov condition for asymptotic stability holds for every detectable pair $(A,C)$. To fill in this gap, next, we use the newly developed notions of \textit{semicontrollability} and \textit{semiobservability} in \cite{HHC:JFI:2011} to give alternative, yet new necessary and sufficient conditions for semistability of $A$ by using the null space condition on $A$ and the Lyapunov equation. Here the definition of semicontrollability is slightly different from the one in \cite{HHC:JFI:2011}.

\begin{definition}
Let $A\in\mathbb{R}^{n\times n}$, $B\in\mathbb{R}^{n\times l}$, and $C\in\mathbb{R}^{l\times n}$.
The pair $(A,B)$ is \textit{semicontrollable} if 
\begin{eqnarray}\label{Range}
{\rm{span}}\left(\bigcup_{i=1}^{n}\mathcal{R}\left(A^{i-1}B\right)\right)=\mathcal{R}(A),
\end{eqnarray} where ${\rm{span}}$ denotes the span of subspace. The pair $(A,C)$ is \textit{semiobservable} if 
\begin{eqnarray}
\bigcap_{i=1}^{n}\mathcal{N}\left(CA^{i-1}\right)=\mathcal{N}(A).\label{Null}
\end{eqnarray}
\end{definition}

As mentioned in \cite{HHC:JFI:2011}, semicontrollability and semiobservability describe the ability of the system to deal with the equilibrium manifold. In particular, semicontrollability implies the ability of the system to drive any initial state to an equilibrium point on the equilibrium manifold in finite time while semiobservability implies that  the only state that is unobservable at the initial time is the equilibrium state. Hence, these two notions are quite different from the classic controllability and observability. However, some properties derived from controllability and observability still hold for these two notions. For example, duality between semicontrollability and semiobservability is valid. Before we present this result, we prove an equivalent test for semicontrollability. 

\begin{lemma}\label{sum_test}
Let $A\in\mathbb{R}^{n\times n}$ and $B\in\mathbb{R}^{n\times l}$. Then $(A,B)$ is semicontrollable if and only if 
\begin{eqnarray}
\sum_{i=1}^{n}\mathcal{R}\left(A^{i-1}B\right)=\mathcal{R}(A).
\end{eqnarray}
\end{lemma}

\begin{IEEEproof}
It suffices to show that $\sum_{i=1}^{n}\mathcal{R}\left(A^{i-1}B\right)={\rm{span}}\left(\bigcup_{i=1}^{n}\mathcal{R}\left(A^{i-1}B\right)\right)$. Note that $\mathcal{R}(A^{i-1}B)$ is a subspace of $\mathbb{R}^{n}$ for every $i=1,\ldots,n$. Then it follows from Fact 2.9.13 in \cite[p.~121]{Bernstein:2009} that the above equality holds.
\end{IEEEproof}

\begin{lemma}\label{subspace}
Let $\mathcal{S}_{1},\mathcal{S}_{2}\subseteq\mathbb{R}^{n}$ be subspaces. Then $\mathcal{S}_{1}=\mathcal{S}_{2}$ if and only if $\mathcal{S}_{1}^{\perp}=\mathcal{S}_{2}^{\perp}$.
\end{lemma}

\begin{IEEEproof}
It follows from $i$) and $ii$) of Fact 2.9.14 in \cite[p.~121]{Bernstein:2009} that $\mathcal{S}_{1}\subseteq\mathcal{S}_{2}$ if and only if $\mathcal{S}_{2}^{\perp}\subseteq\mathcal{S}_{1}^{\perp}$. Hence, the result follows immediately.
\end{IEEEproof}

\begin{lemma}\label{dual}
Let $A\in\mathbb{R}^{n\times n}$ and $C\in\mathbb{R}^{l\times n}$. Then $(A,C)$ is semiobservable if and only if $(A^{\rm{T}},C^{\rm{T}})$ is semicontrollable. 
\end{lemma}

\begin{IEEEproof}
It follows from Theorem 2.4.3 in \cite[p.~103]{Bernstein:2009} that $\mathcal{N}(A)=\mathcal{R}(A^{\rm{T}})^{\perp}$. Hence, if (\ref{Null}) holds, then $\bigcap_{i=1}^{n}\mathcal{R}((A^{\rm{T}})^{i-1}C^{\rm{T}})^{\perp}=\mathcal{R}(A^{\rm{T}})^{\perp}$. By Fact 2.9.16 in \cite[p.~121]{Bernstein:2009}, $\bigcap_{i=1}^{n}\mathcal{R}((A^{\rm{T}})^{i-1}C^{\rm{T}})^{\perp}=[\sum_{i=1}^{n}\mathcal{R}((A^{\rm{T}})^{i-1}\\C^{\rm{T}})]^{\perp}$, and hence, $[\sum_{i=1}^{n}\mathcal{R}((A^{\rm{T}})^{i-1}C^{\rm{T}})]^{\perp}=\mathcal{R}(A^{\rm{T}})^{\perp}$. On the other hand, it follows from Lemma~\ref{subspace} that $\sum_{i=1}^{n}\mathcal{R}((A^{\rm{T}})^{i-1}C^{\rm{T}})=\mathcal{R}(A^{\rm{T}})$. Now it follows from Lemma~\ref{sum_test} that $\sum_{i=1}^{n}\mathcal{R}((A^{\rm{T}})^{i-1}C^{\rm{T}})={\rm{span}}(\bigcup_{i=1}^{n}\mathcal{R}\\((A^{\rm{T}})^{i-1}C^{\rm{T}}))$. This leads to ${\rm{span}}\left(\bigcup_{i=1}^{n}\mathcal{R}\left((A^{\rm{T}})^{i-1}C^{\rm{T}}\right)\right)=\mathcal{R}(A^{\rm{T}})$, and hence, by definition, $(A^{\rm{T}},C^{\rm{T}})$ is semicontrollable. The proof of the converse part just reverses the above derivation. 
\end{IEEEproof}

It is worth noting that semidetectability is weaker than semiobservability.

\begin{lemma}\label{SOSD}
Let $A\in\mathbb{R}^{n\times n}$ and $C\in\mathbb{R}^{l\times n}$. If $(A,C)$ is semiobservable, then $(A,C)$ is semidetectable.  
\end{lemma}

\begin{IEEEproof}
Suppose $(A,C)$ is not semidetectable. Then by definition, there exists $x\neq0$ such that $Cx=0$ and $(A-\jmath\omega I_{n})x=0$ for some $\omega\neq0$. Hence, $CA^{i-1}x=(\jmath\omega)^{i-1}Cx=0$, $i=1,\ldots,n$. Since $(A,C)$ is semiobservable, it follows from (\ref{Null}) that $x\in\mathcal{N}(A)$, and hence, $\jmath\omega x=0$, which leads to $x=0$. This contradicts the assumption. Hence, $(A,C)$ is semidetectable.    
\end{IEEEproof}

The converse part of Lemma~\ref{SOSD} is not true in general as the following example shows. 

\begin{example}
Consider $A$ and $C$ given by Remark~\ref{rmk_sd}. 
It follows from Remark~\ref{rmk_sd} that $(A,C)$ is semidetectable. However, $\mathcal{N}(C)\cap\mathcal{N}(CA)=\mathcal{N}\left(\left[\begin{array}{cc}
0 & 1 \\
\end{array}\right]\right)\neq\mathcal{N}(A)$. Thus, $(A,C)$ is not semiobservable.  \hfill$\blacktriangle$
\end{example}  

Now, combining with Theorem~\ref{thm_rank}, we have the following new series of necessary and sufficient conditions for semistability by using the Lyapunov equation, semiobservability, and semidetectability. 

\begin{theorem}\label{thm_null}
Let $A\in\mathbb{R}^{n\times n}$. Then the following statements are equivalent: 
\begin{itemize}
\item[$i$)] $A$ is semistable.
\item[$ii$)] ${\rm{rank}}\,(A-\jmath\omega I_{n})=n$ for any nonzero $\omega\in\mathbb{R}$ and there exist a positive integer $m$, an $m\times n$ matrix $C$, and an $n\times n$ matrix $P=P^{\rm{T}}>0$ such that (\ref{newLE}) holds.
\item[$iii$)] For every semiobservable pair $(A,C)$ where $C\in\mathbb{R}^{l\times n}$, there exists an $n\times n$
matrix $P=P^{\rm{T}}>0$ such that (\ref{newLE}) holds.
\item[$iv$)] There exist a positive integer $m$, an $m\times n$ matrix $C$, and an $n\times n$ matrix $P=P^{\rm{T}}>0$ such that the pair $(A,C)$ is semiobservable and (\ref{newLE}) holds.
\end{itemize} 
\end{theorem}

\begin{IEEEproof}
First, note that if $A$ is semistable, then it follows from Lemma~\ref{lemma_BP} that $\jmath\omega$, $\omega\neq0$, is not an eigenvalue of $A$. Hence, ${\rm{rank}}\,(A-\jmath\omega I_{n})=n$ for any nonzero $\omega\in\mathbb{R}$. The existence of $C$ and $P$ to (\ref{newLE}) follows from the first part of the proof of Theorem~\ref{thm_rank}. Alternatively, it follows from the second part of the proof of Theorem~\ref{thm_rank} that $A$ is Lyapunov stable. Furthermore, it follows from ${\rm{rank}}\,(A-\jmath\omega I_{n})=n$ for any nonzero $\omega\in\mathbb{R}$ that $\jmath\omega$, $\omega\neq0$, is not an eigenvalue of $A$. Hence, $A$ is semistable. Thus, $i)$ $\Leftrightarrow$ $ii)$.
The proof of $i$) $\Leftrightarrow$ $iii$) follows from Theorem 2.2 in \cite{HHC:JFI:2011}. Next, we show that $i$) $\Leftrightarrow$ $iv$). It follows from Lemma~\ref{SOSD} and Theorem~\ref{thm_rank} that $iv$) $\Rightarrow$ $i$). Alternatively, if $i)$ holds, then we pick $C$ to be the one satisfying $\mathcal{N}(C)=\mathcal{N}(A)$ (an obvious choice is $C=A$). Since $\mathcal{N}(CA^{i})\supseteq\mathcal{N}(A)$ for every $i=1,\ldots,n-1$, it follows that $\bigcap_{i=1}^{n}\mathcal{N}(CA^{i-1})=\mathcal{N}(A)$. Finally, it follows from the similar arguments in the proof of Theorem~\ref{thm_rank} that there exists $P=P^{\rm{T}}>0$ such that (\ref{newLE}) holds. Thus, $i$) $\Rightarrow$ $iv$).    
\end{IEEEproof}

The major difference between Theorems \ref{thm_rank} and \ref{thm_null} is Part $iii$) in Theorem \ref{thm_null}. As we mentioned in Remark~\ref{rmk_sd} and the beginning of this subsection, this part does not hold for semidetectability. Hence, by assuming a stronger notion--semiobservability, one can obtain two necessary and sufficient conditions that are consistent with the classic Lyapunov conditions for asymptotic stability by use of observability.   

\section{Semistable $\mathcal{H}_{2}$ Optimal Control Design for Network Systems}\label{SSH2IAO}

In this section, we use the proposed notions of semistabilizability and semicontrollability in Section~\ref{SSH2} to develop an $\mathcal{H}_{2}$ optimal network-topology-preserving control design methodology for network systems. More specifically, we use Theorems \ref{thm_rank} and \ref{thm_null} to derive the following key result for our design. 

\begin{theorem}\label{main}
Assume that Assumptions \ref{A1} and \ref{A2} hold. 
\begin{itemize}
\item[$i$)] Then $\tilde{A}$ is semistable if and only if for every semicontrollable pair $(\tilde{A},D)$, there exists an $n\times n$ matrix $\tilde{P}=\tilde{P}^{\rm{T}}>-V$ such that 
\begin{eqnarray}\label{LQS_LE}
0=\tilde{A}(\tilde{P}+V)+(\tilde{P}+V)\tilde{A}^{\rm{T}}+DD^{\rm{T}}.
\end{eqnarray}
\item[$ii$)] Assume that there exists an $n\times n$ matrix $\tilde{P}=\tilde{P}^{\rm{T}}>-V$ such that (\ref{LQS_LE}) holds. Then $\tilde{A}$ is semistable if and only if $(\tilde{A},D)$ is semistabilizable. 
\end{itemize}
Furthermore, if $(\tilde{A},D)$ is semistabilizable and $\tilde{P}$ satisfies (\ref{LQS_LE}), then
\begin{eqnarray}\label{hatP}
\tilde{P}=\int_{0}^{\infty}e^{\tilde{A}t}(\tilde{A}V+V\tilde{A}^{\rm{T}}+DD^{\rm{T}})e^{\tilde{A}^{\rm{T}}t}{\rm{d}}t+\alpha xx^{\rm{T}}, \quad x\in\mathcal{N}(\tilde{A}),\quad \alpha>0.
\end{eqnarray}
\end{theorem}

\begin{IEEEproof}
Part $i$) is a direct consequence of Lemma~\ref{dual} and $iii$) of Theorem~\ref{thm_null}.
Part $ii$) is a direct consequence of Theorem~\ref{thm_rank} and $ii$) of Theorem~\ref{thm_null} by letting $P=\tilde{P}+V$. To prove that $\tilde{P}$ has the form (\ref{hatP}), we follow the arguments in the proof of Theorem 4.1 in \cite{Hui:SCL:2011}. More specifically, first it follows from the above result that
there exists an $n\times n$ matrix $\tilde{P}>-V$ such
that (\ref{LQS_LE}) holds or, equivalently, $(\tilde{A}\oplus
\tilde{A}){\rm{vec}}\,\tilde{P}=-{\rm{vec}}\,(\tilde{A}V+V\tilde{A}^{\rm{T}}+DD^{\rm{T}})$. Hence,
${\rm{vec}}\,(\tilde{A}V+V\tilde{A}^{\rm{T}}+DD^{\rm{T}})\in\mathcal{R}(\tilde{A}\oplus\tilde{A})$. Next,
it follows from Lemma 3.8 of \cite{Hui:SCL:2011} that $\tilde{A}\oplus\tilde{A}$ is
semistable, and hence, by Lemma 3.9 of \cite{Hui:SCL:2011},
\begin{eqnarray}
{\rm{vec}}^{-1}\Big((\tilde{A}\oplus
\tilde{A})^{\#}{\rm{vec}}\,(\tilde{A}V+V\tilde{A}^{\rm{T}}+DD^{\rm{T}})\Big)&=&-\int_{0}^{\infty}{\rm{vec}}^{-1}\Big(e^{(\tilde{A}\oplus
\tilde{A})t}{\rm{vec}}\,(\tilde{A}V+V\tilde{A}^{\rm{T}}+DD^{\rm{T}})\Big){\rm{d}}t\nonumber\\
&=&-\int_{0}^{\infty}{\rm{vec}}^{-1}\left(e^{\tilde{A}t}\otimes
e^{\tilde{A}t}\right){\rm{vec}}\,(\tilde{A}V+V\tilde{A}^{\rm{T}}+DD^{\rm{T}}){\rm{d}}t\nonumber\\
&=&-\int_{0}^{\infty}e^{\tilde{A}t}(\tilde{A}V+V\tilde{A}^{\rm{T}}+DD^{\rm{T}})e^{\tilde{A}^{\rm{T}}t}{\rm{d}}t,\label{vec}
\end{eqnarray}
where in (\ref{vec}) we used the facts that $e^{X\oplus
Y}=e^{X}\otimes e^{Y}$ and ${\rm{vec}}(XYZ)=(Z^{\rm{T}}\otimes
X){\rm{vec}}\,Y$ \cite{Bernstein:2009}. Hence,
$\tilde{P}=\int_{0}^{\infty}e^{\tilde{A}t}(\tilde{A}V+V\tilde{A}^{\rm{T}}+DD^{\rm{T}})e^{\tilde{A}t}{\rm{d}}t+{\rm{vec}}^{-1}(z)$,
where $z$ satisfies $z\in\mathcal{N}(\tilde{A}\oplus\tilde{A})$ and
${\rm{vec}}^{-1}(z)=({\rm{vec}}^{-1}(z))^{\rm{T}}\geq 0$ (The nonnegativity of ${\rm{vec}}^{-1}(z)$ is guaranteed by Theorem 4.2a of \cite{SZ:SIAP:1970}). Since $\tilde{A}\oplus\tilde{A}$ is semistable, it follows that the general solution to the equation $(\tilde{A}\otimes\tilde{A})z=0$ is given by $z=x\otimes y$, where $x,y\in\mathcal{N}(\tilde{A})$. Hence, ${\rm{vec}}^{-1}(z)={\rm{vec}}^{-1}(x\otimes y)=yx^{\rm{T}}$, where we used the fact that $xy^{\rm{T}}={\rm{vec}}^{-1}(y\otimes x)$ \cite{Bernstein:2009}. Furthermore, it follows from ${\rm{vec}}^{-1}(z)=({\rm{vec}}^{-1}(z))^{\rm{T}}\geq 0$ that $xy^{\rm{T}}=yx^{\rm{T}}\geq0$. Now it follows from Lemma~\ref{lemma_xy} that $y=\alpha x$, where $\alpha>0$.
\end{IEEEproof}

Using Theorem~\ref{main}, we can obtain a solution to the semistable $\mathcal{H}_{2}$ optimal network-topology-preserving control design proposed in Section~\ref{PFNS} over \textit{some} admissible set by solving the following constrained optimization problem. 

\begin{corollary}\label{thm_2}
Assume that Assumptions \ref{A1} and \ref{A2} hold. Let $\mathcal{K}=\{K: (\tilde{A},D)\,\,\rm{is}\,\,\rm{semistabilizable}\}$. Then solving the following constrained optimization problem
\begin{eqnarray}\label{LMI}
\min_{K\in\mathcal{K}}\Big\{{\rm{tr}}\,(S+V)\tilde{R}:
S=S^{\rm{T}}>-V\,\,{\rm{and}}\,\,\tilde{A}(S+V)+(S+V)\tilde{A}^{\rm{T}}+DD^{\rm{T}}=0\Big\}
\end{eqnarray}  gives a solution of the semistable $\mathcal{H}_{2}$ optimal network-topology-preserving control problem proposed in Section~\ref{PFNS} over $\mathcal{K}$. Furthermore, let $\mathcal{C}_{\rm{s}}=\{K: (\tilde{A},D)\,\,\rm{is}\,\,\rm{semicontrollable}\}$. Then solving the following constrained optimization problem
\begin{eqnarray}\label{LMI_1}
\min_{K\in\mathcal{C}_{\rm{s}}}\Big\{{\rm{tr}}\,(S+V)\tilde{R}:
S=S^{\rm{T}}>-V\,\,{\rm{and}}\,\,\tilde{A}(S+V)+(S+V)\tilde{A}^{\rm{T}}+DD^{\rm{T}}=0\Big\}
\end{eqnarray}  gives a solution of the semistable $\mathcal{H}_{2}$ optimal network-topology-preserving control problem proposed in Section~\ref{PFNS} over $\mathcal{C}_{\rm{s}}$.
\end{corollary}

\begin{IEEEproof}
It follows from Theorem~\ref{main} that $S$ can be parameterized as
$S=\int_{0}^{\infty}e^{\tilde{A}t}(\tilde{A}V+V\tilde{A}^{\rm{T}}+DD^{\rm{T}})e^{\tilde{A}^{\rm{T}}t}{\rm{d}}t+\alpha xx^{\rm{T}}$ 
for some $\alpha>0$ and $x\in\mathcal{N}(\tilde{A})$. Then $S=W+\alpha xx^{\rm{T}}$, which implies that $S+V=W+V+\alpha xx^{\rm{T}}$, where $W$ is given by (\ref{Wmatrix}). Hence,
${\rm{tr}}\,(S+V)\tilde{R}={\rm{tr}}\,(W+V)\tilde{R}+\alpha x^{\rm{T}}\tilde{R}x\geq{\rm{tr}}\,(W+V)\tilde{R}$, where the equality holds if and only if $x\in\mathcal{N}(\tilde{R})\cap\mathcal{N}(\tilde{A})$. The second conclusion follows in a similar fashion.
\end{IEEEproof}

Note that Corollary~\ref{thm_2} is a general optimization framework to solve the robust and optimal control of network systems. The cost functional in (\ref{LMI}) is an implicit function of $K$ due to the fact that $S=S(K)$ satisfies a Lyapunov equation for the closed-loop system.  Hence, it is not easy to solve this constrained optimization problem practically. The authors in \cite{HHC:JFI:2011,HH:IJC:2009} suggest a Linear Matrix Inequality (LMI) method to solve semistable $\mathcal{H}_{2}$ optimal control problems for general linear systems and a Bilinear Matrix Inequality (BMI) method to solve special network consensus problems. However, in general the semistabilizability constraint for $(\tilde{A},D)$ is \textit{not} convex with respect to $K$, and hence, the constrained optimization problem (\ref{LMI}) is \textit{not} a convex optimization problem. Moreover, the existence of a solution to (\ref{LMI}) remains unclear at this point. Thus, it is highly questionable that one can find an optimal solution for (\ref{LMI}) by using the LMI method proposed in \cite{HHC:JFI:2011}.

\section{Optimality Analysis}

To overcome this puzzle of solving (\ref{LMI}), next we use the notion of semicontrollability to show that the constrained optimization problem (\ref{LMI}) indeed has a solution over a subset of $\mathcal{K}$. See also \cite{Haddad:TR} for the discussion of the classic $\mathcal{H}_{2}$ control case. 

\begin{lemma}\label{nonempty}
There always exists $K$ given in the form of (\ref{control_G}) such that $(\tilde{A},D)$ is semicontrollable.
\end{lemma}

\begin{IEEEproof}
We first claim that there exists $K=K(\mathcal{G})$ such that $\mathcal{N}(D^{\rm{T}})=\mathcal{N}(\tilde{A}^{\rm{T}})$. Note that 
\begin{eqnarray}
\tilde{A}_{(i,j)}=\left\{\begin{array}{cc}\mathcal{A}_{(i,j)}(a_{ji}+k_{ji}), & i\neq j \\
\mathcal{A}_{(i,i)}(a_{ii}+k_{ii})-\sum_{r=1,r\neq i}^{n}\mathcal{A}_{(i,r)}(a_{ir}+k_{ir}), & i=j\end{array}\right.
\end{eqnarray} and 
\begin{eqnarray}
D_{(i,j)}=\left\{\begin{array}{cc}\mathcal{A}_{(i,j)}d_{ji}, & i\neq j \\
\mathcal{A}_{(i,i)}d_{ii}-\sum_{r=1,r\neq i}^{n}\mathcal{A}_{(i,r)}d_{ir}, & i=j\end{array}\right..
\end{eqnarray} Clearly if we pick $K$ to be the form $k_{ji}=-a_{ji}+c_{j}d_{ji}$, where $c_{j}\neq0$, $c_{j}\in\mathbb{R}$, $i,j=1,\ldots,n$, then it follows that $\tilde{A}^{\rm{T}}=CD^{\rm{T}}$, where $C={\rm{diag}}[c_{1},\ldots,c_{n}]\in\mathbb{R}^{n\times n}$ is an invertible diagonal matrix. In this case, $\mathcal{N}(\tilde{A}^{\rm{T}})=\mathcal{N}(CD^{\rm{T}})=\mathcal{N}(D^{\rm{T}})$. Finally, note that $\mathcal{N}(D^{\rm{T}}(\tilde{A}^{\rm{T}})^{i})\supseteq\mathcal{N}(\tilde{A}^{\rm{T}})$ for every positive integer $i$, it follows that $\bigcap_{i=1}^{n}\mathcal{N}(D^{\rm{T}}(\tilde{A}^{\rm{T}})^{i-1})=\mathcal{N}(\tilde{A}^{\rm{T}})$, which implies that $(\tilde{A}^{\rm{T}},D^{\rm{T}})$ is semiobservable. By duality (Lemma~\ref{dual}), $(\tilde{A},D)$ is semicontrollable. 
\end{IEEEproof}

\begin{lemma}\label{lemma_Jcont}
Assume that Assumption~\ref{A1} holds. If $K\in\mathcal{S}\cap\mathcal{C}_{\rm{s}}\neq\varnothing$, then $J(\cdot)$ is finite and continuous on $\mathcal{S}\cap\mathcal{C}_{\rm{s}}$.
\end{lemma}

\begin{IEEEproof}
By Lemma~\ref{lemma_int}, if $\tilde{A}$ is semistable, then $Q=\lim_{t\to\infty}Q(t)$ exists if and only if $\mathcal{N}(\tilde{A}^{\rm{T}})\subseteq\mathcal{N}(D^{\rm{T}})$, where $Q(t)$ satisfies (\ref{Qt}) in Lemma~\ref{lemma_Qt}. Since $K\in\mathcal{C}_{\rm{s}}$ implies $\mathcal{N}(\tilde{A}^{\rm{T}})\subseteq\mathcal{N}(D^{\rm{T}})$ by Lemma~\ref{dual}, it follows from Lemma~\ref{lemma_Qt} and (\ref{JKQ}) that $J(\cdot)$ is finite. Finally, the continuity of $J(\cdot)$ follows from (\ref{Qt}). 
\end{IEEEproof}

Recall from \cite{BH:1995} that the pair $(A,B)$ is \textit{completely unstabilizable} if there exists an invertible matrix $S\in\mathbb{R}^{n\times n}$ such that 
\begin{eqnarray}
SAS^{-1}=\left[\begin{array}{cc}
A_{1} & A_{12} \\
0 & A_{2} \\
\end{array}\right], \quad SB=\left[\begin{array}{c}
B_{1} \\
0 \\
\end{array}\right],
\end{eqnarray} where $A_{1}$ is asymptotically stable. The following results give several necessary and sufficient conditions to guarantee the complete unstabilizability of $(A,B)$.

\begin{lemma}\cite[p.~816]{Bernstein:2009}\label{lemma_cud}
The following statements are equivalent: 
\begin{itemize}
\item[$i$)] $(A,B)$ is completely unstabilizable.
\item[$ii$)] $\lim_{t\to\infty}e^{At}B=0$.
\item[$iii$)] $\int_{0}^{\infty}e^{At}BB^{\rm{T}}e^{A^{\rm{T}}t}{\rm{d}}t$ exists.
\item[$iv$)] There exists a positive-semidefinite matrix $Q=Q^{\rm{T}}\in\mathbb{R}^{n\times n}$ satisfying $0=AQ+QA^{\rm{T}}+BB^{\rm{T}}$.
\end{itemize}
\end{lemma} 

The following result can be viewed as a partial converse of Lemma \ref{lemma_Jcont}.

\begin{lemma}\label{lemma_KS}
Assume that Assumptions \ref{A1} and \ref{A2} hold. Let $K\in\mathbb{R}^{m\times n}$. If $J(K)<\infty$, $K\in\mathcal{C}_{\mathrm{s}}\neq\varnothing$, and $\tilde{R}>0$, then $K\in\mathcal{S}$. 
\end{lemma}

\begin{IEEEproof}
By (\ref{cost_J}) and Lemma \ref{lemma_Jcont}, $J(K)=\lim_{t\to\infty}{\rm{tr}}\left[e^{\tilde{A}t}Ve^{\tilde{A}^{\rm{T}}t}\tilde{R}+\int_{0}^{t}e^{\tilde{A}s}DD^{\rm{T}}e^{\tilde{A}^{\rm{T}}s}{\rm{d}}s\tilde{R}\right]$. Note that 
\begin{eqnarray}
e^{\tilde{A}t}Ve^{\tilde{A}^{\rm{T}}t}+\int_{0}^{t}e^{\tilde{A}s}DD^{\rm{T}}e^{\tilde{A}^{\rm{T}}s}{\rm{d}}s&=& V+\int_{0}^{t}\frac{\rm{d}}{{\rm{d}}s}\left(e^{\tilde{A}s}Ve^{\tilde{A}^{\rm{T}}s}\right){{\rm{d}}s}+\int_{0}^{t}e^{\tilde{A}s}DD^{\rm{T}}e^{\tilde{A}^{\rm{T}}s}{\rm{d}}s\nonumber\\
&=& V+\int_{0}^{t}e^{\tilde{A}s}\left(\tilde{A}V+V\tilde{A}^{\rm{T}}+DD^{\rm{T}}\right)e^{\tilde{A}^{\rm{T}}s}{\rm{d}}s.
\end{eqnarray} It follows that $J(K)={\rm{tr}}V\tilde{R}+\lim_{t\to\infty}{\rm{tr}}\left[\int_{0}^{t}e^{\tilde{A}s}\left(\tilde{A}V+V\tilde{A}^{\rm{T}}+DD^{\rm{T}}\right)e^{\tilde{A}^{\rm{T}}s}{\rm{d}}s\tilde{R}\right]$. Thus, $J(K)<\infty$ if and only if $\int_{0}^{\infty}e^{\tilde{A}s}\left(\tilde{A}V+V\tilde{A}^{\rm{T}}+DD^{\rm{T}}\right)e^{\tilde{A}^{\rm{T}}s}{\rm{d}}s$ exists. Next, if $P=\int_{0}^{\infty}e^{\tilde{A}s}\left(\tilde{A}V+V\tilde{A}^{\rm{T}}+DD^{\rm{T}}\right)e^{\tilde{A}^{\rm{T}}s}{\rm{d}}s$ exists, then $P=P^{\rm{T}}\geq-V$ and $\lim_{t\to\infty}e^{\tilde{A}t}\left(\tilde{A}V+V\tilde{A}^{\rm{T}}+DD^{\rm{T}}\right)e^{\tilde{A}^{\rm{T}}t}=0$. In this case,  
\begin{eqnarray}\label{JWV1}
\tilde{A}(P+V)+(P+V)\tilde{A}^{\rm{T}}&=&\int_{0}^{\infty}\left[\tilde{A}e^{\tilde{A}s}\left(\tilde{A}V+V\tilde{A}^{\rm{T}}+DD^{\rm{T}}\right)e^{\tilde{A}^{\rm{T}}s}\right.\nonumber\\
&&\left.+e^{\tilde{A}s}\left(\tilde{A}V+V\tilde{A}^{\rm{T}}+DD^{\rm{T}}\right)e^{\tilde{A}^{\rm{T}}s}\tilde{A}^{\rm{T}}\right]{\rm{d}}s\nonumber\\
&&+\tilde{A}V+V\tilde{A}^{\rm{T}}\nonumber\\
&=&\int_{0}^{\infty}\frac{{\rm{d}}}{{\rm{d}}s}e^{\tilde{A}s}\left(\tilde{A}V+V\tilde{A}^{\rm{T}}+DD^{\rm{T}}\right)e^{\tilde{A}^{\rm{T}}s}{\rm{d}}s\nonumber\\
&&+\tilde{A}V+V\tilde{A}^{\rm{T}}\nonumber\\
&=&\lim_{t\to\infty}e^{\tilde{A}t}\left(\tilde{A}V+V\tilde{A}^{\rm{T}}+DD^{\rm{T}}\right)e^{\tilde{A}^{\rm{T}}t}\nonumber\\
&&-(\tilde{A}V+V\tilde{A}^{\rm{T}}+DD^{\rm{T}})+\tilde{A}V+V\tilde{A}^{\rm{T}}\nonumber\\
&=&-DD^{\rm{T}}.
\end{eqnarray} Now it follows from Lemma~\ref{lemma_cud} that the pair $(\tilde{A},D)$ is completely unstabilizable and $\int_{0}^{\infty}e^{\tilde{A}s}DD^{\rm{T}}e^{\tilde{A}^{\rm{T}}s}{\rm{d}}s$ exists.

Consider the pair $(\tilde{A},D)$. Then it follows from the Kalman decomposition that there exists an invertible matrix $T\in\mathbb{R}^{n\times n}$ such that $T\tilde{A}T^{-1}=\left[\begin{array}{cc}
\overline{A}_{1} & \overline{A}_{12} \\
0 & \overline{A}_{2} \\
\end{array}\right]$ and $TD=\left[\begin{array}{c}
D_{1} \\
0 \\
\end{array} \right]$ where $(\overline{A}_{1},D_{1})$ is controllable. Thus,
$\int_{0}^{\infty}e^{\tilde{A}s}DD^{\rm{T}}e^{\tilde{A}^{\rm{T}}s}{\rm{d}}s=T^{-1}\left[\begin{array}{cc}
\int_{0}^{\infty}e^{\overline{A}_{1}s}D_{1}D_{1}^{\rm{T}}e^{\overline{A}_{1}^{\rm{T}}s}{\rm{d}}s & 0 \\ 
0 & 0 
\end{array} \right](T^{-1})^{\rm{T}}$ exists. 

First, we claim that $\overline{A}_{1}$ is asymptotically stable. Conversely, suppose that $\overline{A}_{1}$ is not asymptotically stable. Let $\lambda\in{\rm{spec}}(\tilde{A})$ where ${\rm{Re}}\lambda\geq0$, and let $x\in\mathbb{C}^{r}$, $x\neq0$ satisfy $\overline{A}_{1}x=\lambda x$, where $r$ denotes the row dimension of $\overline{A}_{1}$. Since $(\overline{A}_{1},D_{1})$ is controllable, it follows from Theorem 12.6.18 of \cite[p.~815]{Bernstein:2009} that $\int_{0}^{\infty}e^{\overline{A}_{1}s}D_{1}D_{1}^{\rm{T}}e^{\overline{A}_{1}^{\rm{T}}s}{\rm{d}}s$ is positive definite. Thus, $\alpha\triangleq x^{*}\int_{0}^{\infty}e^{\overline{A}_{1}s}D_{1}D_{1}^{\rm{T}}e^{\overline{A}_{1}^{\rm{T}}s}{\rm{d}}sx$ is a positive number. However, we also have $\alpha=x^{*}\int_{0}^{\infty}e^{\lambda s}D_{1}D_{1}^{\rm{T}}e^{\overline{\lambda} s}{\rm{d}}sx=x^{*}D_{1}D_{1}^{\rm{T}}x\int_{0}^{\infty}e^{2({\rm{Re}}\lambda)s}{\rm{d}}s$. By controllability of $(\overline{A}_{1},D_{1})$, we have $x^{*}D_{1}D_{1}^{\rm{T}}x\neq0$. Since ${\rm{Re}}\lambda\geq0$, it follows that $\int_{0}^{\infty}e^{2({\rm{Re}}\lambda)s}{\rm{d}}s=\infty$, which contradicts the fact that $\alpha$ is a positive number. Thus, $\overline{A}_{1}$ is asymptotically stable.

Second, we claim that $\tilde{A}$ is Lyapunov stable. To see this, let $\mu\in{\rm{spec}}(\tilde{A}^{\mathrm{T}})$ and $z\in\mathbb{C}^{n}$, $z\neq0$ satisfy $\tilde{A}^{\mathrm{T}}z=\mu z$. Since $(\tilde{A},D)$ is semicontrollable, it follows from Theorem~\ref{main} that $z^{*}DD^{\rm{T}}z=-z^{*}(\tilde{A}(\hat{P}+V)+(\hat{P}+V)\tilde{A}^{\mathrm{T}})z=-(2{\mathrm{Re}}\mu)z^{*}(\hat{P}+V)z$ for some $\hat{P}=\hat{P}^{\mathrm{T}}>-V$. Hence, ${\rm{Re}}\mu\leq0$. If ${\rm{Re}}\mu=0$, then let $\mu=\jmath\omega$ and $x\in\mathcal{N}((\jmath\omega I_{n}-\tilde{A}^{\mathrm{T}})^{2})$, where $\omega\in\mathbb{R}$. Defining $y=(\jmath\omega I_{n}-\tilde{A}^{\mathrm{T}})x$, it follows that $(\jmath\omega I_{n}-\tilde{A}^{\mathrm{T}})y=0$, and hence, $\tilde{A}^{\mathrm{T}}y=\jmath\omega y$. Therefore, it follows from Theorem~\ref{main} that $-y^{*}DD^{\rm{T}}y=y^{*}(\tilde{A}(\hat{P}+V)+(\hat{P}+V)\tilde{A}^{\mathrm{T}})y=-\jmath\omega y^{*}(\hat{P}+V)y+\jmath\omega y^{*}(\hat{P}+V)y=0$, and thus, $D^{\rm{T}}y=0$. Hence, $0=x^{*}DD^{\rm{T}}y=-x^{*}(\tilde{A}(\hat{P}+V)+(\hat{P}+V)\tilde{A}^{\mathrm{T}})y=-x^{*}(\tilde{A}+\jmath\omega I_{n})(\hat{P}+V)y=y^{*}(\hat{P}+V)y$. Since $\hat{P}+V>0$, it follows that $y=0$, that is, $(\jmath\omega I_{n}-\tilde{A}^{\mathrm{T}})x=0$. Therefore, $x\in\mathcal{N}(\jmath\omega I_{n}-\tilde{A}^{\mathrm{T}})$. Now it follows from Proposition 5.5.8 of \cite[p.~323]{Bernstein:2009} that $\jmath\omega$ is semisimple.  
Thus, $\tilde{A}$ is Lyapunov stable and consequently, so is $\overline{A}_{2}$.

Finally, we show that $\jmath\omega\not\in{\rm{spec}}(\overline{A}_{2})$ for every $\omega\neq0$. Since $(\tilde{A},D)$ is semistabilizable, it follows from (\ref{dim}) that $\dim[\mathcal{R}(D)+\mathcal{R}(\jmath\omega I_{n}-\tilde{A})]=n$ for every nonzero $\omega\in\mathbb{R}$. We claim that $(T\tilde{A}T^{-1},TD)$ is semistabilizable. Note that $T$ is invertible, it follows that $\mathcal{R}(TD)=\mathcal{R}(D)$ and for every nonzero $\omega\in\mathbb{R}$, $\mathcal{R}(\jmath\omega I_{n}-T\tilde{A}T^{-1})=\mathcal{R}(\jmath\omega I_{n}-\tilde{A})$. Hence,  $\dim[\mathcal{R}(TD)+\mathcal{R}(\jmath\omega I_{n}-T\tilde{A}T^{-1})]=\dim[\mathcal{R}(D)+\mathcal{R}(\jmath\omega I_{n}-\tilde{A})]=n$, which implies that ${\rm{rank}}\left[\begin{array}{cc}
TD & \jmath\omega I_{n}-T\tilde{A}T^{-1}
\end{array}\right]=n$, that is, ${\rm{rank}}\left[\begin{array}{ccc}
D_{1} & \jmath\omega I-\overline{A}_{1} & -\overline{A}_{12} \\ 
0 & 0 & \jmath\omega I-\overline{A}_{2}
\end{array} \right]=n$. In this case, the matrix $\jmath\omega I-\overline{A}_{2}$ must be a full rank matrix for every nonzero $\omega$, which implies that $\jmath\omega$ is not an eigenvalue of $\overline{A}_{2}$.  Hence, the eigenvalue of $\overline{A}_{2}$ is either a real/complex number with negative real part or is zero, and if the eigenvalue is zero, then it is semisimple. Therefore, $\overline{A}_{2}$ is semistable. Since $\tilde{A}$ has the same set of eigenvalues as $\overline{A}_{1}$ and $\overline{A}_{2}$, it follows that $\tilde{A}$ is semistable, i.e., $K\in\mathcal{S}$. 
\end{IEEEproof}

To proceed, let $\overline{\sigma}(S)$ denote the largest singular value of $S$, $\underline{\sigma}(S)$ the smallest singular value of $S$, $\overline{\lambda}(S)$ the most positive eigenvalue of symmetric $S=S^{\rm{T}}$, and $\underline{\lambda}(S)$ the most negative eigenvalue of symmetric $S=S^{\rm{T}}$.

\begin{lemma}\label{lemma_AVQ}
Suppose that there exist $A\in\mathbb{R}^{n\times n}$, $V=V^{\rm{T}}\in\mathbb{R}^{n\times n}$, and $Q=Q^{\rm{T}}\in\mathbb{R}^{n\times n}$ satisfying $0=AQ+QA^{\rm{T}}+V$. Then $\underline{\lambda}(V)\leq2\underline{\sigma}(Q)\overline{\sigma}(A)$.
\end{lemma}

\begin{IEEEproof}
For every $\alpha>0$, since $0\leq(\sqrt{\alpha}I_{n}+\frac{1}{\sqrt{\alpha}}AQ)(\sqrt{\alpha}I_{n}+\frac{1}{\sqrt{\alpha}}AQ)^{\rm{T}}=\alpha I_{n}+AQ+QA^{\rm{T}}+\frac{1}{\alpha}AQ^{2}A^{\rm{T}}$, it follows that $V=-AQ-QA^{\rm{T}}\leq\alpha I_{n}+\frac{1}{\alpha}AQ^{2}A^{\rm{T}}$. Hence, it follows from Theorem 8.4.9 of \cite[p.~469]{Bernstein:2009} that $\underline{\lambda}(V)\leq\underline{\lambda}(\alpha I_{n}+\frac{1}{\alpha}AQ^{2}A^{\rm{T}})=\alpha+\frac{1}{\alpha}\underline{\lambda}(AQ^{2}A^{\rm{T}})=\alpha+\frac{1}{\alpha}\underline{\sigma}^{2}(AQ)$. Now using the fact that $\underline{\sigma}(AB)\leq\underline{\sigma}(A)\overline{\sigma}(B)$, we have $\underline{\lambda}(V)\leq\alpha+\frac{1}{\alpha}\underline{\sigma}^{2}(Q)\overline{\sigma}^{2}(A)$. Suppose either $\underline{\sigma}(Q)=0$ or $\overline{\sigma}(A)=0$. Then it follows that $\underline{\lambda}(V)\leq\alpha$ for any $\alpha>0$. Hence, $\underline{\lambda}(V)\leq0$, as required in this case. If, however, $\underline{\sigma}(Q)>0$ and $\overline{\sigma}(A)>0$, then set $\alpha=\underline{\sigma}(Q)\overline{\sigma}(A)>0$. Hence, $\underline{\lambda}(V)\leq\underline{\sigma}(Q)\overline{\sigma}(A)+\frac{1}{\underline{\sigma}(Q)\overline{\sigma}(A)}\underline{\sigma}^{2}(Q)\overline{\sigma}^{2}(A)=2\underline{\sigma}(Q)\overline{\sigma}(A)$, as required. 
\end{IEEEproof}

For the next result define the set 
\begin{eqnarray}
\mathcal{K}_{\alpha}\triangleq\{K\in\mathcal{C}_{\rm{s}}: J(K)\leq\alpha\}
\end{eqnarray} and note that, by Lemma~\ref{lemma_KS}, $\mathcal{K}_{\alpha}\subset\mathcal{S}$, $\alpha\geq0$, provided that $\tilde{R}>0$.

\begin{lemma}\label{lemma_compact}
Assume that Assumptions \ref{A1} and \ref{A2} hold. Furthermore, assume $R_{1}>0$ and $\mathcal{C}_{\rm{s}}\neq\varnothing$. Then there exists $\alpha>0$ such that $\mathcal{K}_{\alpha}$ is nonempty and compact relative to $\mathcal{C}_{\rm{s}}$.
\end{lemma} 

\begin{IEEEproof}
First note that it follows from Lemma~\ref{nonempty} that $\mathcal{C}_{\rm{s}}$ is nonempty and from Lemma~\ref{SOSD} that $\mathcal{C}_{\rm{s}}\subseteq\mathcal{K}$. Thus, there exists $\alpha>0$ such that $\mathcal{K}_{\alpha}$ is nonempty.
Define the function $\hat{J}:\mathcal{C}_{\rm{s}}\to\mathbb{R}$ by 
\begin{eqnarray}
\hat{J}(K)\triangleq\left\{\begin{array}{cc}
J(K), & K\in\mathcal{K}_{2\alpha}, \\ 
2\alpha, & K\not\in\mathcal{K}_{2\alpha}.
\end{array} \right.
\end{eqnarray} Since, by Lemma~\ref{lemma_KS}, $\mathcal{K}_{2\alpha}\subset\mathcal{S}$, it follows from Lemma~\ref{lemma_Jcont} that $J(\cdot)$ is continuous on $\mathcal{K}_{2\alpha}$. However, $\hat{J}(K)=J(K)$ for all $K\in\mathcal{K}_{2\alpha}$. Then it follows that $\hat{J}(\cdot)$ is continuous on $\mathcal{K}_{2\alpha}$. Next note that $J(K)\to2\alpha$ as $K\to\partial{\mathcal{K}}_{2\alpha}$. Hence, $\hat{J}(\cdot)$ is continuous on $\mathcal{C}_{\rm{s}}$. Thus $\hat{\mathcal{K}}_{\alpha}\triangleq\{K\in\mathcal{C}_{\rm{s}}:\hat{J}(K)\leq\alpha\}$ is closed relative to $\mathcal{C}_{\rm{s}}$. On the other hand, $\hat{J}(K)\leq\alpha$ implies that $\hat{J}(K)=J(K)$. Then it follows that $\hat{\mathcal{K}}_{\alpha}=\mathcal{K}_{\alpha}$ and thus $\mathcal{K}_{\alpha}$ is closed relative to $\mathcal{C}_{\rm{s}}$. Next we show that $\mathcal{K}_{\alpha}$ is bounded relative to $\mathcal{C}_{\rm{s}}$. Since $\mathcal{K}_{\alpha}\subset\mathcal{S}$, it follows from (\ref{JWV}) that $K\in\mathcal{K}_{\alpha}$ implies that $J(K)={\rm{tr}}(W+V)\tilde{R}$, where $W+V\geq0$ and $0=\tilde{A}(W+V)+(W+V)\tilde{A}^{\rm{T}}+DD^{\rm{T}}$. Now, since ${\rm{tr}}(W+V)\tilde{R}\leq\alpha$ for $K\in\mathcal{K}_{\alpha}$, it follows from $\underline{\lambda}(W+V)=\underline{\sigma}(W+V)$,  $\underline{\lambda}(R_{2})=\underline{\sigma}(R_{2})$, and Proposition 8.4.13 of \cite[p.~471]{Bernstein:2009} that  
$\alpha\geq{\rm{tr}}(W+V)\tilde{R}\geq({\rm{tr}}\tilde{R})\underline{\sigma}(W+V)=({\rm{tr}}R_{1}+{\rm{tr}}K^{\rm{T}}R_{2}K)\underline{\sigma}(W+V)=({\rm{tr}}R_{1}+{\rm{tr}}KK^{\rm{T}}R_{2})\underline{\sigma}(W+V)\geq({\rm{tr}}R_{1}+({\rm{tr}}KK^{\rm{T}})\underline{\sigma}(R_{2}))\underline{\sigma}(W+V)$. By using Lemma~\ref{lemma_AVQ}, we have $({\rm{tr}}R_{1}+({\rm{tr}}KK^{\rm{T}})\underline{\sigma}(R_{2}))\underline{\sigma}(W+V)\geq({\rm{tr}}R_{1}+\|K\|_{\rm{F}}^{2}\underline{\sigma}(R_{2}))\underline{\sigma}(DD^{\rm{T}})/(2\overline{\sigma}(\tilde{A}))$, where $\|K\|_{\rm{F}}$ denotes the Frobenius norm of $K$. Note that $\overline{\sigma}(A+BK)\leq\overline{\sigma}(A)+\overline{\sigma}(B)\overline{\sigma}(K)$ and $\overline{\sigma}(K)\leq\|K\|_{\rm{F}}$. Then it follows that $({\rm{tr}}R_{1}+\|K\|_{\rm{F}}^{2}\underline{\sigma}(R_{2}))\underline{\sigma}(DD^{\rm{T}})/(2\overline{\sigma}(\tilde{A}))\geq([\underline{\sigma}(R_{2})\|K\|_{\rm{F}}^{2}+{\rm{tr}}R_{1}]\underline{\sigma}(DD^{\rm{T}}))/(2[\overline{\sigma}(A)+\overline{\sigma}(B)\overline{\sigma}(K)])\geq([\underline{\sigma}(R_{2})\|K\|_{\rm{F}}^{2}+{\rm{tr}}R_{1}]\underline{\sigma}(V))/(2[\overline{\sigma}(B)\|K\|_{\rm{F}}+\overline{\sigma}(A)])=(\underline{\sigma}(DD^{\rm{T}})\underline{\sigma}(R_{2})\|K\|_{\rm{F}}^{2}+\underline{\sigma}(DD^{\rm{T}}){\rm{tr}}R_{1})/(2\overline{\sigma}(B)\|K\|_{\rm{F}}+2\overline{\sigma}(A))$. Thus, combining these inequalities together yields $\underline{\sigma}(DD^{\rm{T}})\underline{\sigma}(R_{2})\|K\|_{\rm{F}}^{2}+\underline{\sigma}(DD^{\rm{T}}){\rm{tr}}R_{1}\leq2\alpha\overline{\sigma}(B)\|K\|_{\rm{F}}+2\alpha\overline{\sigma}(A)$. Define $\beta_{1}\triangleq(2\alpha\overline{\sigma}(B))/(\underline{\sigma}(DD^{\rm{T}})\underline{\sigma}(R_{2}))$ and $\beta_{2}\triangleq(\underline{\sigma}(DD^{\rm{T}}){\rm{tr}}R_{1}-2\alpha\overline{\sigma}(A))/(\underline{\sigma}(DD^{\rm{T}})\underline{\sigma}(R_{2}))$. Choose $\alpha>0$ to be sufficiently large so that $\beta_{1}^{2}\geq4\beta_{2}$. Then it follows that $\|K\|_{\rm{F}}\leq\frac{1}{2}\beta_{1}+\frac{1}{2}\sqrt{\beta_{1}^{2}-4\beta_{2}}<\infty$, which establishes the boundedness of $\mathcal{K}_{\alpha}$ relative to $\mathcal{C}_{\rm{s}}$. Hence, $\mathcal{K}_{\alpha}$ is compact relative to $\mathcal{C}_{\rm{s}}$.
\end{IEEEproof}

\begin{theorem}\label{thm_existence}
Assume that Assumptions \ref{A1} and \ref{A2} hold. Furthermore, assume $R_{1}>0$ and $\mathcal{S}\cap\mathcal{C}_{\rm{s}}\neq\varnothing$. Then there exists $K_{*}\in\mathcal{S}\cap\mathcal{C}_{\rm{s}}$ such that $J(K_{*})\leq J(K)$, $K\in\mathcal{C}_{\rm{s}}$. 
\end{theorem}

\begin{IEEEproof} 
Since $\mathcal{S}\cap\mathcal{C}_{\rm{s}}\neq\varnothing$, let $\hat{K}\in\mathcal{S}\cap\mathcal{C}_{\rm{s}}\neq\varnothing$ and define $\hat{\alpha}\triangleq J(\hat{K})$. Then it follows from $\hat{K}\in\mathcal{K}_{\hat{\alpha}}$ that $\mathcal{K}_{\hat{\alpha}}\neq\varnothing$ and, by Lemma~\ref{lemma_compact}, $\mathcal{K}_{\hat{\alpha}}$ is compact relative to $\mathcal{C}_{\rm{s}}$. Since $\mathcal{K}_{\hat{\alpha}}\subset\mathcal{S}\cap\mathcal{C}_{\rm{s}}$ by Lemma~\ref{lemma_KS}, it follows from Lemma~\ref{lemma_Jcont} that $J(\cdot)$ is continuous on $\mathcal{K}_{\hat{\alpha}}$. Hence, there exists $K_{*}\in\mathcal{K}_{\hat{\alpha}}$ such that $J(K_{*})\leq J(K)$, $K\in\mathcal{K}_{\hat{\alpha}}$. Furthermore, since $K\not\in\mathcal{K}_{\hat{\alpha}}$ implies that $J(K)>\hat{\alpha}$, it follows that $J(K_{*})\leq J(K)$, $K\in\mathcal{C}_{\rm{s}}$. 
\end{IEEEproof}

Since it has been proved in Theorem~\ref{thm_2} that solving the semistable $\mathcal{H}_{2}$ optimal network-topology-preserving control problem proposed in Section~\ref{PFNS} is equivalent to solving (\ref{LMI}) over the admissible set $\mathcal{C}_{\mathrm{s}}$, it follows from Theorem~\ref{thm_existence} that the following constrained optimization problem
\begin{eqnarray}
\min_{K\in\mathcal{C}_{\rm{s}}}\Big\{{\rm{tr}}\,(S+V)\tilde{R}:
S=S^{\rm{T}}>-V\,\,{\rm{and}}\,\,\tilde{A}(S+V)+(S+V)\tilde{A}^{\rm{T}}+DD^{\rm{T}}=0\Big\}
\end{eqnarray} has an optimal solution and gives a solution to the original semistable $\mathcal{H}_{2}$ optimal network-topology-preserving control problem proposed in Section~\ref{PFNS} over $\mathcal{C}_{\mathrm{s}}$. Thus, finding an appropriate numerical algorithm to solve (\ref{LMI}) becomes a key issue to tackle the proposed optimal and robust control design for network systems.

\section{Numerical Algorithms}

In this section, we propose a heuristic swarm optimization based numerical algorithm to solve (\ref{LMI}). In particular, we use an variant form of Particle Swarm Optimization (PSO) proposed by Kennedy and Eberhart \cite{KE:ICNN:1995} to deal with constrained optimization problems. The standard particle swarm
optimization algorithm in vector form is described as follows:
\begin{eqnarray}
\textbf{v}_{k}(t+1)&=& a\textbf{v}_{k}(t)+b_{1}r_{1}(\textbf{p}_{1,k}-\textbf{x}_{k}(t))+b_{2}r_{2}(\textbf{p}_{2}-\textbf{x}_{k}(t)),\nonumber\\
\textbf{x}_{k}(t+1)&=&\textbf{x}_{k}(t)+\textbf{v}_{k}(t+1),
\end{eqnarray} where
$\textbf{v}_{k}(t)$ and $\textbf{x}_{k}(t)$ are the velocity and
position of particle $k$ at iteration $t$ respectively,
$\textbf{p}_{1,k}$ is the position of previous best value particle
$k$ obtained so far, $\textbf{p}_{2}$ is the position of the global
best value the swarm can achieve so far, $a$, $b_{1}$, and $b_{2}$
are weight coefficients, and $r_{1}$, $r_{2}$ are two random
coefficients which are usually selected in uniform distribution in
the set $\Omega\subseteq[0,1]$. In every iteration $t$, the velocity of each
particle is updated by the interaction of the current velocity, the
previous best position $\textbf{p}_{1,k}$ and the global position
$\textbf{p}_{2}$. The position of each particle is updated by using
the current position and the newly updated velocity. For every
iteration, the previous best position $\textbf{p}_{1,k}$ and the
global best position $\textbf{p}_{2}$ will be updated according to
the calculation of the objective function.

The constrained optimization problem that we are interested in can
be formulated as follows: Find $\textbf{x}$ which minimizes
$f(\textbf{x})$ subject to $g_i(\textbf{x})\leq 0$, $i=1,\ldots,m$, and 
$h_j(\textbf{x})=0$, $j=1,\ldots,p$, 
where $\textbf{x}\in \mathbb{R}^n$ is the solution vector and each
$x_i$, $i=1,\ldots,n$, is bounded by the lower and upper limits
$L_i\leq x_i\leq U_i$ which define the search space, $\mathscr{F}$
comprises the set of all solutions which satisfy the constraints of
the problems and it is called the feasible region, $m$ is the number
of inequality constraints, and $p$ is the number of equality
constraints. Equality constraints are always transformed into
inequality constrains in practice as follows:
$|h_j(\textbf{x})|-\varepsilon \leq 0$ \cite{probdef,LRMCSCD:TR}, where
$\varepsilon$ is the tolerance allowed (a very small value).

It is well known that PSO is designed for unconstrained optimization problems. Regarding the application of PSO to constrained optimization,
there are also a number of papers in the literature
\cite{bestpso,LRMCSCD:TR,bilevel}. In these constrained PSO algorithms, the updates of
$\textbf{p}_{1,k}$ and $\textbf{p}_{2}$ are different from those in
unconstrained PSO. More specifically, three rules to update the
global best position $\textbf{p}_{2}$ and previous best position
$\textbf{p}_{1,k}$ are as follows \cite{3rules}:
\begin{itemize}
\item[1.] If $\textbf{x}_1$ and $\textbf{x}_2$ are both feasible solutions,
then choose the one which can obtain the best objective function
value.
\item[2.] If $\textbf{x}_1$ and $\textbf{x}_2$ are not feasible, then choose the
one with the lowest sum of normalized constraint violation:
\begin{eqnarray}
\hspace{-2em}s(\textbf{x})=\sum^m_{i=1}\max(0,g(\textbf{x}))+\sum^p_{j=1}\max(0,
(|h(\textbf{x})|-\varepsilon)).
\end{eqnarray}
\item[3.] If one is feasible and the other is infeasible, then choose the
feasible one.
\end{itemize}

In this report, since $V$ is an arbitrary symmetric matrix related to initial conditions, without loss of generality, we assume that
$V=0$. Because the pair $(\tilde{A},D)$ is semistabilizable if and
only if $\mbox{rank}[D,\jmath\omega I_{n}-\tilde{A}]=n$, we can use
the rank value to check the semistabilizability of $(\tilde{A},D)$. Furthermore, we use the indicator $\mathfrak{p}=\mathfrak{p}(S)$ of the command $\mathtt{chol}(S)=[R,\mathfrak{p}]$ (Cholesky factorization) in MATLAB$^{\textregistered}$ to test the positive
definiteness of matrix $S$, that is, if $\mathfrak{p}=0$ then $S$ is
positive definite, and if $\mathfrak{p}>0$, then $S$ is not
positive definite.

Then we can formulate our problem into the standard constrained
optimization problem as follows:
\begin{eqnarray}
&&\min f(K)=\min\mbox{tr}(S\tilde{R})\\
&&\hspace{-4em}\mbox{subject to} \nonumber\\
&& h_{1}(K)=\mbox{rank}[D,\jmath\omega I_{n}-\tilde{A}]-n=0,\\
&& h_{2}(K)=\mathfrak{p}=0,\\
&& h_{3}(K)=\tilde{A}S+S\tilde{A}^{\mathrm{T}}+DD^{\mathrm{T}}=0.
\end{eqnarray}

Since it is not easy to make the solution feasible and the matrix
equation makes the calculation complicated in this problem, we
propose a hierarchical method to solve this problem. We divide the
optimization into two stages and in each stage we use a constrained
PSO algorithm to solve the problem. In this report, we also use the those
``feasibility-rules" mentioned above to handle constraints. The proposed numerical algorithm is described in Algorithm~\ref{NGPSO}.

\begin{algorithm}[ht]
\caption{}\label{NGPSO}
\begin{center}
\begin{algorithmic}
  \STATE Initialize matrix particles $K(\mathcal{G})$ in the search
  space with $k_{ij}$ being a uniform
  distribution $U(-\theta,\theta)$, $\theta>0$.\\
\textbf{STAGE I}: \\
  \REPEAT
  \STATE Use the constrained PSO to optimize the problem
  with single constraint $h_{1}(K)=0$.
  \UNTIL{$h_{1}(K)=0$ or the exit condition is satisfied.}
  \STATE Solve the matrix equation $h_{3}(K)=0$ for $S$ and go to STAGE II.\\
\textbf{STAGE II}:\\
  \REPEAT
   \STATE Use the constrained PSO with constraint $h_{2}(K)=0$.
If $h_{1}(K)\neq0$ then go back to STAGE I
  \UNTIL{All the particles converge to the same position or the exit condition is satisfied.}
\end{algorithmic}
\end{center}
\end{algorithm}

\section{Simulation Results}

\subsection{Case 1: $\Omega=[0,1]$}

In this case, we choose a four-node graph topology, a
six-node graph topology, and a ten-node graph topology to test algorithms and give the
numerical results respectively. Before we proceed, we need to define three new weight matrices as
$\mathscr{A}$ and $\mathscr{D}$ in (\ref{ADK}) where each item is
the weights $a_{ij}$ and $d_{ij}$ in $A(\mathcal{G})$ and $D(\mathcal{G})$, respectively. For the 4-node network, we choose 30 particles and
iterate 100 times. For the 6-node network we pick 100 particles and
iterate 100 times. Finally, for the 10-node network, the number of particles
is set to 2000 and we choose 2000 iteration times. The reason that
we pick our particles and iteration times in such a way is that as
the number of nodes increases it is more difficult to find a
feasible solution. Finally, let $\theta=100$ in Algorithm~\ref{NGPSO}.
\begin{figure*}
\begin{eqnarray}\label{ADK}
\mathscr{A}=
\begin{bmatrix}
a_{11} & a_{12} & \cdots & a_{1n}\\
a_{21} & a_{22} & \cdots & a_{2n}\\
\vdots & \vdots & \ddots & \vdots\\
a_{n1} & a_{n2} & \cdots & a_{nn}
\end{bmatrix},\mathscr{D}=
\begin{bmatrix}
d_{11} & d_{12} & \cdots & d_{1n}\\
d_{21} & d_{22} & \cdots & d_{2n}\\
\vdots & \vdots & \ddots & \vdots\\
d_{n1} & d_{n2} & \cdots & d_{nn}
\end{bmatrix}
\end{eqnarray}
\end{figure*}

For the node-4 graph, the topology of the graph is shown in Fig. \ref{fig1}, and $\mathcal{A}$, $\mathscr{A}$, $\mathscr{D}$, $E_{1}$, and $E_{2}$ are given by (\ref{4ADK}).
\begin{figure*}
\begin{eqnarray}\label{4ADK}
\mathcal{A}=
\begin{bmatrix}
\begin{smallmatrix}
1 & 1 & 0 & 1\\
1 & 1 & 0 & 1\\
0 & 0 & 1 & 1\\
1 & 1 & 1 & 1
\end{smallmatrix}
\end{bmatrix},\mathscr{A}=
\begin{bmatrix}
\begin{smallmatrix}
1 & 3 & 0 & 2\\
3 & 1 & 0 & 1\\
0 & 0 & 1 & 2\\
2 & 1 & 2 & 1
\end{smallmatrix}
\end{bmatrix},\mathscr{D}=
\begin{bmatrix}
\begin{smallmatrix}
1 & 1 & 0 & 2\\
1 & 1 & 0 & 3\\
0 & 0 & 1 & 2\\
2 & 3 & 2 & 1
\end{smallmatrix}
\end{bmatrix},E_{1}=
\begin{bmatrix}
\begin{smallmatrix}
1 & 3 & 5 & 7\\
-1 & -3 & -5 & -7
\end{smallmatrix}
\end{bmatrix},E_{2}=
\begin{bmatrix}
\begin{smallmatrix}
2 & 4 & 6 & 8\\
2 & 4 & 6 & 8
\end{smallmatrix}
\end{bmatrix}
\end{eqnarray}
\end{figure*}
\begin{figure}
\centering
\includegraphics[scale=0.4]{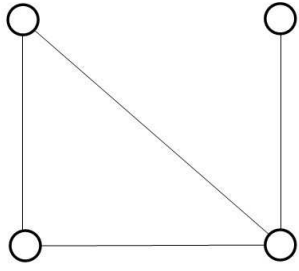}
\caption{The 4-node graph topology}\label{fig1}
\end{figure}
For the node-6 graph, the topology of the graph is shown in Fig. \ref{fig2}, and $\mathcal{A}$, $\mathscr{A}$, $\mathscr{D}$, $E_{1}$, and $E_{2}$ are given by (\ref{6ADK}).
\begin{figure}
\centering
\includegraphics[scale=0.4]{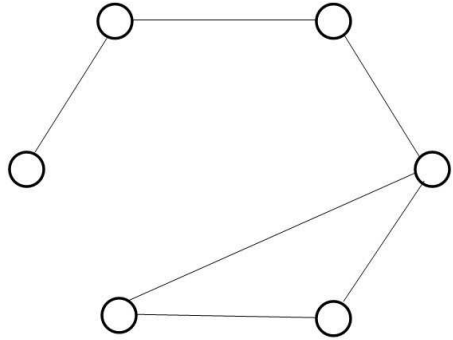}
\caption{The 6-node graph topology}\label{fig2}
\end{figure}
\begin{figure*}
\begin{eqnarray}\label{6ADK}
\mathcal{A}=
\begin{bmatrix}
\begin{smallmatrix}
1 & 1 & 0 & 0 & 0 & 0\\
1 & 1 & 1 & 0 & 0 & 0\\
0 & 1 & 1 & 1 & 0 & 0\\
0 & 0 & 1 & 1 & 1 & 1\\
0 & 0 & 0 & 1 & 1 & 1\\
0 & 0 & 0 & 1 & 1 & 1
\end{smallmatrix}
\end{bmatrix},\mathscr{A}=
\begin{bmatrix}
\begin{smallmatrix}
1 & 1 & 0 & 0 & 0 & 0\\
1 & 1 & 1 & 0 & 0 & 0\\
0 & 1 & 1 & 1 & 0 & 0\\
0 & 0 & 1 & 1 & 1 & 1\\
0 & 0 & 0 & 1 & 1 & 1\\
0 & 0 & 0 & 1 & 1 & 1
\end{smallmatrix}
\end{bmatrix},\mathscr{D}=
\begin{bmatrix}
\begin{smallmatrix}
1 & 1 & 0 & 0 & 0 & 0\\
1 & 1 & 1 & 0 & 0 & 0\\
0 & 1 & 1 & 1 & 0 & 0\\
0 & 0 & 1 & 1 & 1 & 1\\
0 & 0 & 0 & 1 & 1 & 1\\
0 & 0 & 0 & 1 & 1 & 1
\end{smallmatrix}
\end{bmatrix},E_{1}=
\begin{bmatrix}
\begin{smallmatrix}
1 & 3 & 5 & 7 & 9 & 11\\
-1 & -3 & -5 & -7 & -9 & -11
\end{smallmatrix}
\end{bmatrix},E_{2}=
\begin{bmatrix}
\begin{smallmatrix}
2 & 4 & 6 & 8 & 10 & 12\\
2 & 4 & 6 & 8 & 10 & 12
\end{smallmatrix}
\end{bmatrix}
\end{eqnarray}
\end{figure*}
For the node-10 graph, the topology of the graph is shown in Fig. \ref{fig3}, and $\mathcal{A}$, $\mathscr{A}$, $\mathscr{D}$, $E_{1}$, and $E_{2}$ are given by (\ref{10ADK1}) and (\ref{10ADK2}). 
\begin{figure}
\centering
\includegraphics[scale=0.4]{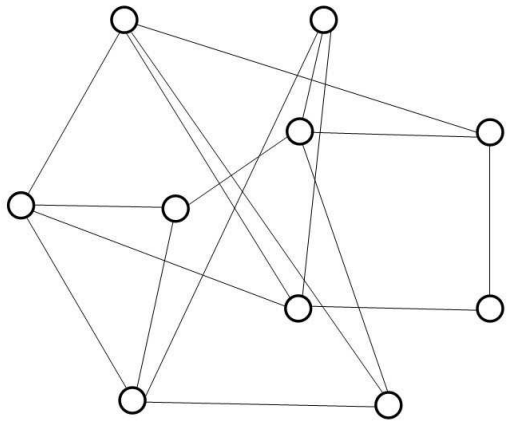}
\caption{The 10-node graph topology}\label{fig3}
\end{figure}
\begin{figure*}
\begin{eqnarray}
\mathcal{A}=
\begin{bmatrix}
\begin{smallmatrix}
1 & 1 & 0 & 1 & 0 & 0 & 1 & 0 & 1 & 0 \\
1&1&1&0&0&0&0&1&1&0\\
0&1&1&0&1&0&1&1&0&0\\
1&0&0&1&0&0&0&1&0&1\\
0&0&1&0&1&1&0&0&0&1\\
0&0&0&0&1&1&0&1&0&0\\
1&0&1&0&0&0&1&0&0&1\\
0&1&1&1&0&1&0&1&0&0\\
1&1&0&0&0&0&0&0&1&1\\
0&0&0&1&1&0&1&0&1&1
\end{smallmatrix}
\end{bmatrix},\mathscr{A}=
\begin{bmatrix}
\begin{smallmatrix}
1&4&0&0&0&0&3&0&4&0\\
4&1&5&0&0&0&0&3&4&0\\
0&5&1&0&3&0&1&2&0&0\\
0&0&0&1&0&0&0&3&0&4\\
0&0&3&0&1&4&0&0&0&5\\
0&0&0&0&4&1&0&3&0&0\\
3&0&1&0&0&0&1&0&0&3\\
0&3&2&3&0&3&0&1&0&0\\
4&4&0&0&0&0&0&0&1&4\\
0&0&0&4&5&0&3&0&4&1
\end{smallmatrix}
\end{bmatrix},\mathscr{D}=
\begin{bmatrix}
\begin{smallmatrix}
1&2&0&0&0&0&3&0&2&0\\
2&1&1&0&0&0&0&3&2&0\\
0&1&1&0&3&0&5&4&0&0\\
0&0&0&1&0&0&0&3&0&2\\
0&0&3&0&1&2&0&0&0&1\\
0&0&0&0&2&1&0&3&0&0\\
3&0&5&0&0&0&1&0&0&3\\
0&3&4&3&0&3&0&1&0&0\\
2&2&0&0&0&0&0&0&1&4\\
0&0&0&2&1&0&3&0&4&1
\end{smallmatrix}
\end{bmatrix}\label{10ADK1}\\
E_{1}=
\begin{bmatrix}
\begin{smallmatrix}
1&3&5&7&9&11&13&15&17&19\\
-1&-3&-5&-7&-9&-11&-13&-15&-17&-19
\end{smallmatrix}
\end{bmatrix},E_{2}=
\begin{bmatrix}
\begin{smallmatrix}
2 & 4 & 6 & 8 & 10 & 12 & 14 & 16 & 18 & 20\\
2 & 4 & 6 & 8 & 10 & 12 & 14 & 16 & 18 & 20
\end{smallmatrix}
\end{bmatrix}\label{10ADK2}
\end{eqnarray}
\end{figure*}

The numerical results are listed as follows: 

1) For the 4-node graph,
the best value we found is $4.1499\times 10^{4}$ and the solution
is given by (\ref{4K}).
\begin{figure*}
\begin{eqnarray}\label{4K}
K=
\begin{bmatrix}
-25.8128&6.9356&0&7.4131\\
-15.2475&-48.7597&0&28.8212\\
0&0&16.9044&39.1779\\
33.7429&9.9462&-39.9050&-63.8953\\
\end{bmatrix}
\end{eqnarray} 
\end{figure*}
We set the initial value $x_{0}=\mbox{[$-1$ $-5$ $3$ $2$]}^{\mathrm{T}}$
to the system and simulate the response of the system in Fig. \ref{fig4}.
\begin{figure}
\centering
\includegraphics[scale=0.55]{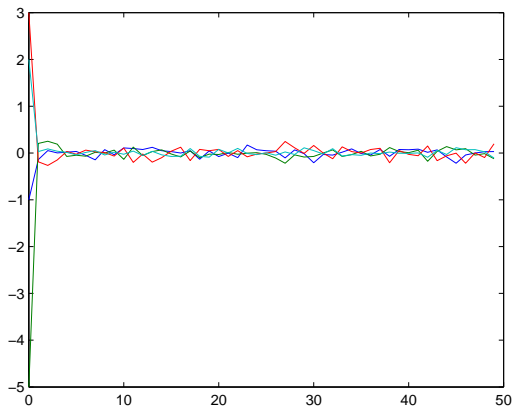}
\caption{The 4-node network system}\label{fig4}
\end{figure}

2) For the 6-node graph,
the best value we found is $6.5001\times 10^{4}$ and the solution
is given by (\ref{6K}).
\begin{figure*}
\begin{eqnarray}\label{6K}
K=
\begin{bmatrix}
-5.7999&63.3691&0&0&0&0\\
-47.9432&-5.3760&-12.3541&0&0&0\\
0&-11.6398&-19.2400&-20.1679&0&0\\
0&0&89.0270&-260.7268&71.1626&32.1778\\
0&0&0&36.6739&-89.8199&3.2043\\
0&0&0&64.7487&-95.3711&-73.1023
\end{bmatrix}
\end{eqnarray}
\end{figure*}
We set the initial value $x_{0}=\mbox{[$-1$ $2$ $-3$ $-5$ $6$
$4$]}^{\mathrm{T}}$ to the system and simulate the response of the
system in Fig. \ref{fig5}.
\begin{figure}
\centering
\includegraphics[scale=0.55]{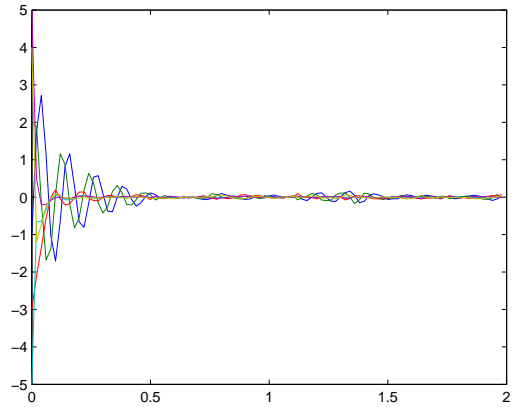}
\caption{The 6-node network system}\label{fig5}
\end{figure}

3) For the 10-node graph,
the best value we found is $1.6542\times 10^{7}$ and the solution
is given by (\ref{10K}).
\begin{figure*}
\begin{eqnarray}\label{10K}
K=
\begin{bmatrix}
\begin{smallmatrix}
-150.5268&5.9712&0&55.7382&0&0&47.5353&0&16.4830&0\\
7.5120&-40.5108&59.9242&0&0&0&0&34.1820&-64.0819&0\\
0&56.9838&-127.6148&0&-26.2837&0&47.0278&54.8586&0&0\\
-57.5445&0&0&-84.8205&0&0&0&48.6216&0&48.0671\\
0&0&24.7723&0&-106.9229&-5.2833&0&0&0&83.4365\\
0&0&0&0&91.7352&-84.7957&0&-5.3865&0&0\\
-100.6877&0&87.7404&0&0&0&-26.6525&0&0&-51.2242\\
0&-26.3584&27.2982&39.1605&0&43.4890&0&-185.2530&0&0\\
-70.3768&4.2835&0&0&0&0&0&0&-16.6587&64.6209\\
0&0&0&9.2189&59.6637&0&38.6097&0&-29.5420&-94.3604
\end{smallmatrix}
\end{bmatrix}
\end{eqnarray}
\end{figure*}
We set the initial value $x_{0}=\mbox{[$5$ $-5$ $3$ $-1$ $1$ $2$ $-3$
$4$ $-2$ $-6$]}^{\mathrm{T}}$ to the system and simulate the
response of the system in Fig. \ref{fig6}. 
\begin{figure}
\centering
\includegraphics[scale=0.55]{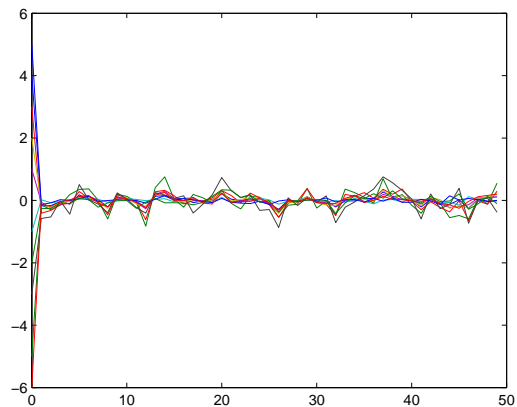}
\caption{The 10-node network system}\label{fig6}
\end{figure}

\subsection{Case 2: $\Omega=\{0.01k,k=1,\ldots,100\}$}

In this case, we first choose a 2-node network system with $x_{i}(t)\in\mathbb{R}^{2}$, $i=1,2$, and then we apply the proposed algorithm to two
additional network systems with particular topologies given by \cite{OFM:PIEEE:2007} and \cite{Strogatz:Nature:2001}, respectively, for obtaining
numerical results. Let $\Omega=\{0.01k,k=1,\ldots,100\}$. For the 2-node network with the 2-dimensional state, we choose 200 particles and
iterate 100 times. For the 20-node and 24-node networks with 1-dimensional states, the number of particles
is set to 2000 and we choose 2000 iteration times. The reason that
we pick our particles and iteration times in such a way is that as
the number of nodes increases it is more difficult to find a
feasible solution. 

\begin{figure}
\centering
\epsfig{file=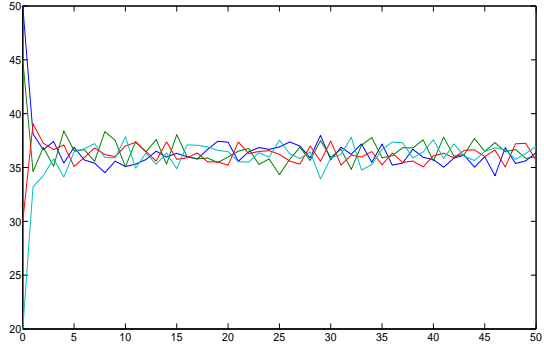,width=0.5\linewidth,height=0.3\linewidth}
\caption{State versus time for the 2-node network system with the two-dimensional state for each node.}\label{fig2agents}
\end{figure}

In the 2-node network system with $x_{i}(t)\in\mathbb{R}^{2}$,
$\mathcal{A}=
\begin{bmatrix}
\begin{smallmatrix}
0&1\\
1&0
\end{smallmatrix}
\end{bmatrix}$, $a_{ii}=d_{ii}=0$, $a_{12}=-a_{21}=-1$, $d_{12}=-d_{21}=-1$, $i,j=1,2$, 
$E_{1}=\begin{bmatrix}
\begin{smallmatrix}
1&0&0&1\\
1&2&2&1\\
1&0&1&4\\
1&1&1&2\\
1&0&0&1
\end{smallmatrix}
\end{bmatrix}$, and $E_{2}=\begin{bmatrix}
\begin{smallmatrix}
1&1&1&1\\
0&0&0&0\\
0&0&0&0\\
0&0&0&0\\
-1&-1&-1&-1
\end{smallmatrix}
\end{bmatrix}$.
The best value we found for $f(K)$ is $84.625$.
and the
corresponding matrix $K$ is given by
\begin{eqnarray*}
K=
\begin{bmatrix}
\begin{smallmatrix}
-94.1723&34.7539&35.2913&24.1272\\
34.7539&-122.8293  & 98.1980 & -10.1225\\
35.2913&98.198&-150.9285&17.4393\\
24.1272 & -10.1225 &  17.4393 & -31.4439 \end{smallmatrix}
\end{bmatrix}.
\end{eqnarray*}
The time response of the closed-loop system is given by Figure \ref{fig2agents}. Note that in this case, $x_{\infty}\neq0$.

\begin{figure}
\hfill
\begin{minipage}{.49\linewidth}
\centering
\includegraphics[height=0.6\textwidth, width=0.7\textwidth]{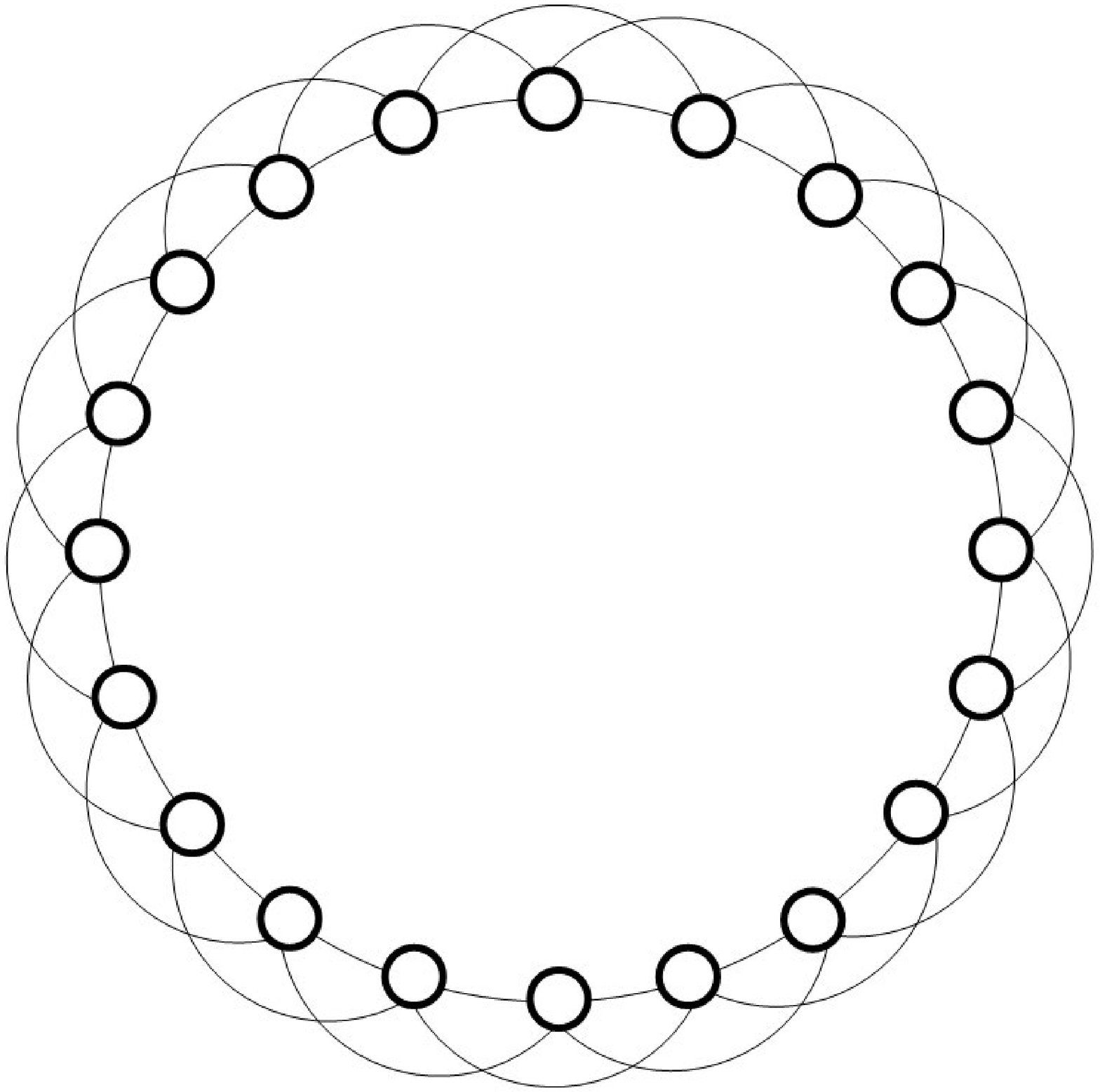} \\
\caption{The 20-node graph topology.}\label{fig20nodes}
\end{minipage}
\hfill
\begin{minipage}{.49\linewidth}
\centering
\includegraphics[height=0.6\textwidth, width=0.7\textwidth]{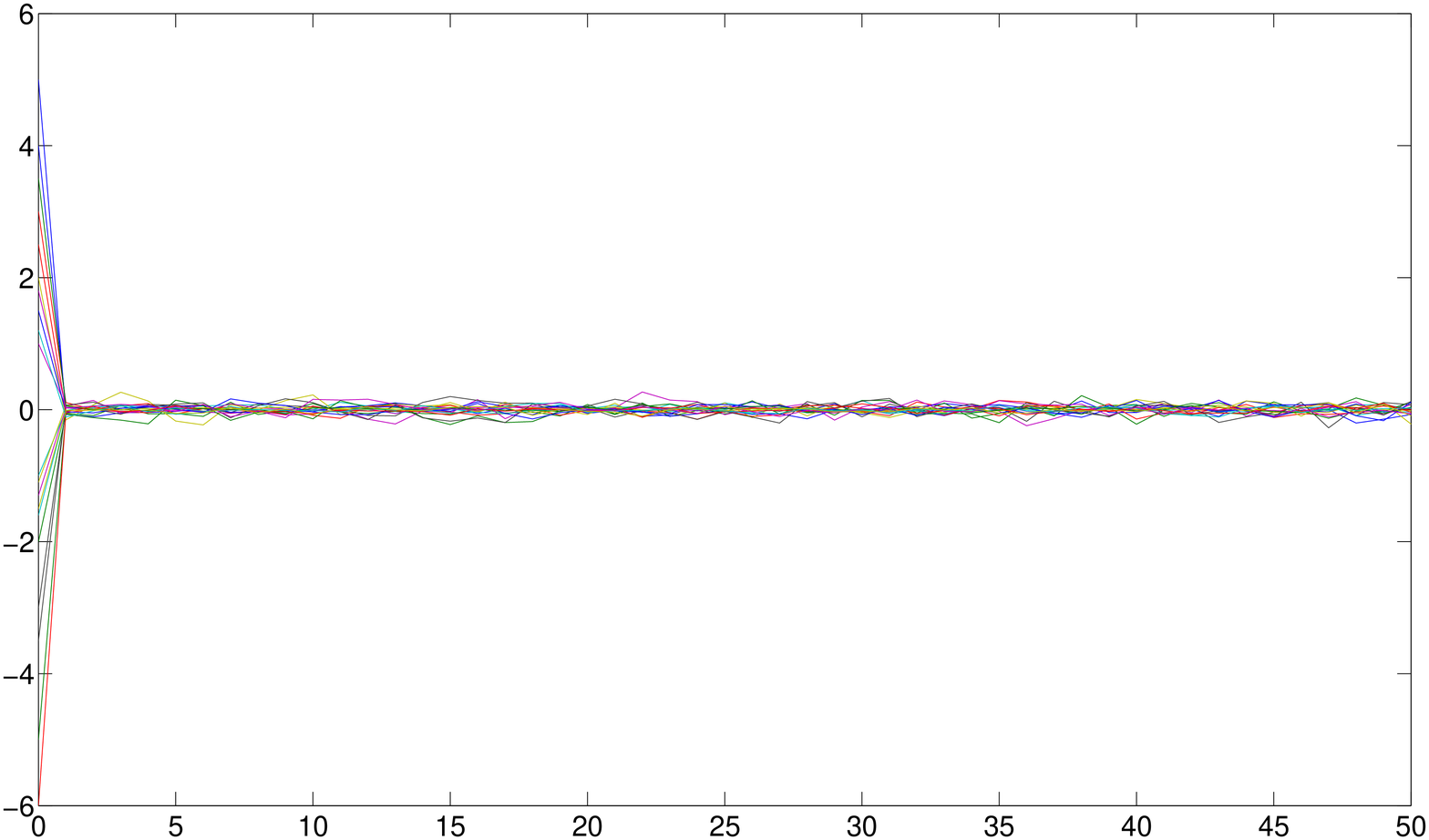} \\
\caption{State versus time of the 20-node network system.}\label{figrps20node}
\end{minipage}
\hfill
\end{figure}

\begin{figure}
\hfill
\begin{minipage}{.49\linewidth}
\centering
\includegraphics[height=0.6\textwidth, width=0.7\textwidth]{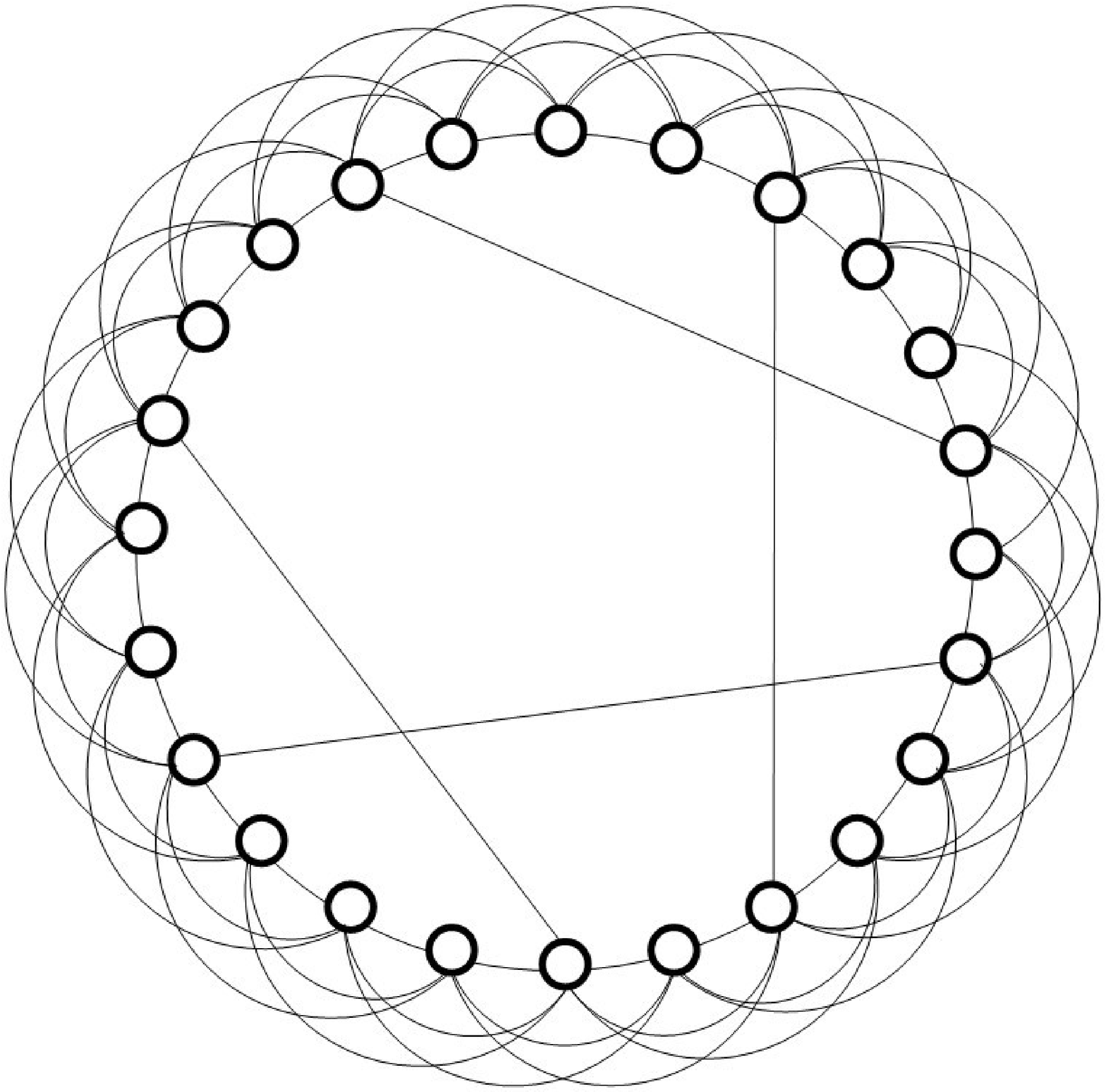} \\
\caption{The 24-node graph topology.}\label{fig24nodes}
\end{minipage}
\hfill
\begin{minipage}{.49\linewidth}
\centering
\includegraphics[height=0.6\textwidth, width=0.7\textwidth]{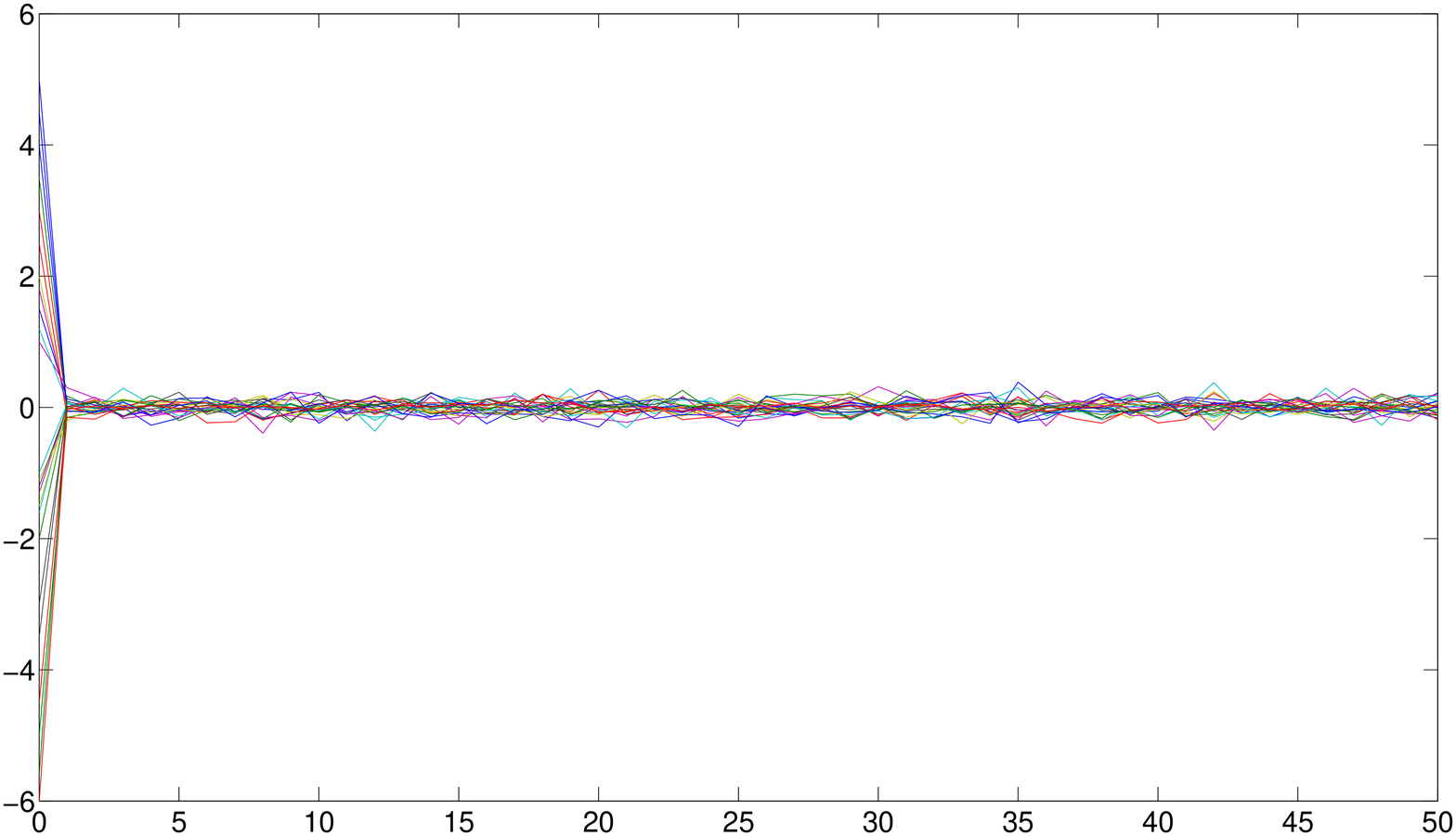} \\
\caption{State versus time of the 24-node network system.}\label{figrps24node}
\end{minipage}
\hfill
\end{figure}

For the 20-node and 24-node network systems, the topologies are
shown in Figures \ref{fig20nodes} and \ref{fig24nodes}, which are originally from \cite{OFM:PIEEE:2007} and \cite{Strogatz:Nature:2001}, respectively.
Due to the space limitation, we do not present the corresponding
matrices $R_{1}$ and $R_{2}$ as well as $a_{ij}$ and $d_{ij}$ we used in the simulation here. For the 20-node graph, the best value we found for $f(K)$ is $1906.2$ and the simulation of the state versus time for this network system is shown in
Figure \ref{figrps20node}. Finally, for the 24-node graph, the best value we found for $f(K)$ is $2772$ and the simulation of the state versus time for this network system is shown in
Figure \ref{figrps24node}. Note that for these two network systems, $x_{\infty}=0$.

\section{Conclusion}

Motivated by multiagent network consensus and deterministic optimal semistable control, this report lays out
a semistable $\mathcal{H}_{2}$ control synthesis framework called LQS for stochastic linear network systems via the
notions of semistabilizability and semicontrollability. By exploiting the properties of semistability and
the relationship between semidetectability and the Lyapunov equation, a complete set of necessary and sufficient conditions were
established for semistability and a constrained optimization characterization of semistable $\mathcal{H}_{2}$ control over some admissible set was derived. Next, we have addressed the existence of optimal solutions to such a constrained optimization problem. Furthermore, we developed a constrained PSO based numerical algorithm to solve the proposed optimization problem. The future work focuses on the application of the proposed numerical method to design a robust and optimal consensus protocol with imperfect information for multiagent networks proposed in \cite{Hui:ACC:2012}. 

\appendix

\begin{lemma}\label{limit}
If $f:(-\infty,+\infty)\to\mathbb{R}$ is integrable and $\lim_{t\to+\infty}f(t)=A$ where $-\infty\leq A\leq+\infty$, then
\begin{eqnarray*}
\lim_{t\to+\infty}\frac{1}{t}\int_{0}^{t}f(s){\mathrm{d}}s=A=\lim_{t\to+\infty}f(t).
\end{eqnarray*}
\end{lemma}

\begin{IEEEproof}
We first consider the case where $-\infty<A<+\infty$. It follows from $\lim_{t\to+\infty}f(t)=A$ that for every $\varepsilon>0$, there exists $\delta>0$ such that for all $t>\delta$, $|f(t)-A|<\varepsilon$, or, equivalently, $-\varepsilon+A<f(t)<\varepsilon+A$. Hence, for all $t>\delta$, $\frac{1}{t}\int_{0}^{\delta}f(s){\mathrm{d}}s+\frac{t-\delta}{t}(-\varepsilon+A)<\frac{1}{t}\int_{0}^{\delta}f(s){\mathrm{d}}s+\frac{1}{t}\int_{\delta}^{t}f(s){\mathrm{d}}s<(\varepsilon+A)\frac{t-\delta}{t}+\frac{1}{t}\int_{0}^{\delta}f(s){\mathrm{d}}s$, that is, $\frac{1}{t}\int_{0}^{\delta}f(s){\mathrm{d}}s-\frac{t-\delta}{t}\varepsilon-\frac{\delta}{t}A<\frac{1}{t}\int_{0}^{t}f(s){\mathrm{d}}s-A<\frac{t-\delta}{t}\varepsilon-\frac{\delta}{t}A+\frac{1}{t}\int_{0}^{\delta}f(s){\mathrm{d}}s$ for all $t>\delta$, and hence, $\frac{1}{t}\int_{0}^{\delta}f(s){\mathrm{d}}s-\varepsilon-\frac{\delta}{t}A<\frac{1}{t}\int_{0}^{t}f(s){\mathrm{d}}s-A<\varepsilon-\frac{\delta}{t}A+\frac{1}{t}\int_{0}^{\delta}f(s){\mathrm{d}}s$ for all $t>\delta$. Note that $\lim_{t\to+\infty}(-\frac{\delta}{t}A+\frac{1}{t}\int_{0}^{\delta}f(s){\mathrm{d}}s)=0$, it follows that there exists $\hat{\delta}>0$ such that for all $t>\hat{\delta}$, $|-\frac{\delta}{t}A+\frac{1}{t}\int_{0}^{\delta}f(s){\mathrm{d}}s|<\varepsilon$, or, equivalently, $-\varepsilon<-\frac{\delta}{t}A+\frac{1}{t}\int_{0}^{\delta}f(s){\mathrm{d}}s<\varepsilon$. Thus, for $t>\max\{\delta,\hat{\delta}\}$, $-2\varepsilon<\frac{1}{t}\int_{0}^{t}f(s){\mathrm{d}}s-A<2\varepsilon$. By definition, $\lim_{t\to+\infty}\frac{1}{t}\int_{0}^{t}f(s){\mathrm{d}}s=A$.

Now we consider the case where $A=+\infty$. Conversely, suppose that $\lim_{t\to+\infty}\frac{1}{t}\int_{0}^{t}f(s){\mathrm{d}}s\neq+\infty$. Since $\lim_{t\to+\infty}f(t)=+\infty$, it follows that $\int_{0}^{t}f(s){\mathrm{d}}s$ must have a lower bound for all $t\geq0$, and hence, $\liminf_{t\to+\infty}\\\frac{1}{t}\int_{0}^{t}f(s){\mathrm{d}}s\geq0$, i.e., there exists $\delta_{1}>0$ such that $\frac{1}{t}\int_{0}^{t}f(s){\mathrm{d}}s\geq0$ for all $t>\delta_{1}$. Consequently, if $\lim_{t\to+\infty}\frac{1}{t}\int_{0}^{t}f(s){\mathrm{d}}s<+\infty$, then there exists $M\in[0,+\infty)$ such that $\lim_{t\to+\infty}\frac{1}{t}\int_{0}^{t}f(s){\mathrm{d}}s=M$, which implies that $\limsup_{t\to+\infty}\frac{1}{t}\int_{0}^{t}f(s){\mathrm{d}}s=M$. On the other hand, it follows from $\lim_{t\to+\infty}f(t)=+\infty$ that there exists $\delta_{2}>0$ such that $f(t)>2M+1$ for all $t>\delta_{2}$. Let $\delta=\max\{\delta_{1},\delta_{2}\}$. Then for every $t>2\delta+1$,  $\frac{1}{t}\int_{0}^{t}f(s){\mathrm{d}}s=\frac{1}{t}\int_{0}^{\delta}f(s){\mathrm{d}}s+\frac{1}{t}\int_{\delta}^{t}f(s){\mathrm{d}}s=\frac{\delta}{t}\frac{1}{\delta}\int_{0}^{\delta}f(s){\mathrm{d}}s+\frac{1}{t}\int_{\delta}^{t}f(s){\mathrm{d}}s\geq\frac{1}{t}\int_{\delta}^{t}f(s){\mathrm{d}}s\geq\frac{t-\delta}{t}(2M+1)>\frac{1}{2}(2M+1)=M+\frac{1}{2}$, which implies that $\limsup_{t\to+\infty}\frac{1}{t}\int_{0}^{t}f(s){\mathrm{d}}s>M$. This contradicts the fact that $\limsup_{t\to+\infty}\frac{1}{t}\int_{0}^{t}f(s){\mathrm{d}}s=M$. Thus, $\lim_{t\to+\infty}\frac{1}{t}\int_{0}^{t}f(s){\mathrm{d}}s=+\infty$. 

Finally, for the case where $A=-\infty$, define $g(t)=-f(t)$. Then it follows from $\lim_{t\to+\infty}f(t)=-\infty$ that $\lim_{t\to+\infty}g(t)=+\infty$. By the second case, $\lim_{t\to+\infty}\frac{1}{t}\int_{0}^{t}g(s){\mathrm{d}}s=+\infty$, which implies that $-\lim_{t\to+\infty}\frac{1}{t}\int_{0}^{t}\\f(s){\mathrm{d}}s=+\infty$. Hence, $\lim_{t\to+\infty}\frac{1}{t}\int_{0}^{t}f(s){\mathrm{d}}s=-\infty$.
\end{IEEEproof}

\begin{lemma}\label{lemma_xy}
Let $x,y\in\mathbb{R}^{n}$ be two nonzero vectors. Then $xy^{\mathrm{T}}=yx^{\mathrm{T}}$ if and only if $x$ and $y$ are linearly dependent. Furthermore, $xy^{\mathrm{T}}=yx^{\mathrm{T}}\geq0$ if and only if $y=\alpha x$, where $\alpha>0$.
\end{lemma}

\begin{IEEEproof}
If $x$ and $y$ are linearly dependent, then $xy^{\mathrm{T}}=yx^{\mathrm{T}}$ holds. Conversely, suppose that $x$ and $y$ are linearly independent, then it follows from Proposition 7.1.8 of \cite[p.~441]{Bernstein:2009} that $xy^{\mathrm{T}}=yx^{\mathrm{T}}$ if and only if ${\rm{vec}}^{-1}(y\otimes x)={\rm{vec}}^{-1}(x\otimes y)$, which further implies that $y\otimes x=x\otimes y$, where ${\mathrm{vec}}^{-1}$ denotes the inverse operation of vectorization \cite[p.~439]{Bernstein:2009}. Let $x=[x_{1},\ldots,x_{n}]^{\mathrm{T}}$ and $y=[y_{1},\ldots,y_{n}]^{\mathrm{T}}$. Then it follows from $y\otimes x=x\otimes y$ that $x_{i}x=x_{i}y$ for every $i=1,\ldots,n$. Since $x$ and $y$ are linearly independent, it follows that $y_{i}x-x_{i}y=0$ for every $i=1,\ldots,n$ if and only if $y_{i}=x_{i}=0$ for every $i=1,\ldots,n$. This contradicts the assumption that $x,y\neq 0$. Hence, $x$ and $y$ are linearly dependent. Similar arguments hold for the second result.
\end{IEEEproof}


\bibliographystyle{IEEEtran}
\bibliography{Reference}

\begin{thebibliography}{10}
\providecommand{\url}[1]{#1}
\csname url@samestyle\endcsname
\providecommand{\newblock}{\relax}
\providecommand{\bibinfo}[2]{#2}
\providecommand{\BIBentrySTDinterwordspacing}{\spaceskip=0pt\relax}
\providecommand{\BIBentryALTinterwordstretchfactor}{4}
\providecommand{\BIBentryALTinterwordspacing}{\spaceskip=\fontdimen2\font plus
\BIBentryALTinterwordstretchfactor\fontdimen3\font minus
  \fontdimen4\font\relax}
\providecommand{\BIBforeignlanguage}[2]{{%
\expandafter\ifx\csname l@#1\endcsname\relax
\typeout{** WARNING: IEEEtran.bst: No hyphenation pattern has been}%
\typeout{** loaded for the language `#1'. Using the pattern for}%
\typeout{** the default language instead.}%
\else
\language=\csname l@#1\endcsname
\fi
#2}}
\providecommand{\BIBdecl}{\relax}
\BIBdecl

\bibitem{HH:IJC:2009}
Q.~Hui and W.~M. Haddad, ``$\mathcal{H}_{2}$ optimal semistable stabilization
  for linear discrete-time dynamical systems with applications to network
  consensus,'' \emph{Int. J. Control}, vol.~82, pp. 456--469, 2009.

\bibitem{Hui:SCL:2011}
Q.~Hui, ``Optimal semistable control for continuous-time linear systems,''
  \emph{Syst. Control Lett.}, vol.~60, no.~4, pp. 278--284, 2011.

\bibitem{HHC:JFI:2011}
W.~M. Haddad, Q.~Hui, and V.~Chellaboina, ``$\mathcal{H}_{2}$ optimal
  semistable control for linear dynamical systems: {A}n {LMI} approach,''
  \emph{J. Franklin Inst.}, vol. 348, no.~10, pp. 2898--2910, 2011.

\bibitem{Hui:AUT:2011}
Q.~Hui, ``Optimal distributed linear averaging,'' \emph{Automatica}, vol.~47,
  no.~12, pp. 2713--2719, 2011.

\bibitem{Hui:JFI:2012}
------, ``Distributed semistable {LQR} control for discrete-time dynamically
  coupled systems,'' \emph{J. Franklin Inst.}, vol. 349, no.~1, pp. 74--92,
  2012.

\bibitem{HHB:TAC:2008}
Q.~Hui, W.~M. Haddad, and S.~P. Bhat, ``Finite-time semistability and consensus
  for nonlinear dynamical networks,'' \emph{IEEE Trans. Autom. Control},
  vol.~53, pp. 1887--1900, 2008.

\bibitem{BB:JVA:1995}
D.~S. Bernstein and S.~P. Bhat, ``Lyapunov stability, semistability, and
  asymptotic stability of matrix second-order systems,'' \emph{ASME J. Vibr.
  Acoustics}, vol. 117, pp. 145--153, 1995.

\bibitem{HHB:ACC:2008}
Q.~Hui, W.~M. Haddad, and S.~P. Bhat, ``Semistability theory for differential
  inclusions with applications to consensus problems in dynamical networks with
  switching topology,'' in \emph{Proc. Amer. Control Conf.}, Seattle, WA, 2008,
  pp. 3981--3986.

\bibitem{Hui:MTNS:2010}
Q.~Hui, ``Optimal finite-time distributed linear averaging,'' in \emph{19th
  Int. Symp. Math. Theory Networks Syst.}, Budapest, Hungary, 2010, pp.
  1961--1968.

\bibitem{Hui:TAC:2011}
------, ``Finite-time rendezvous algorithms for mobile autonomous agents,''
  \emph{IEEE Trans. Autom. Control}, vol.~56, no.~1, pp. 207--211, 2011.

\bibitem{SDD:CDC:2007}
H.~Sandberg, J.-C. Delvenne, and J.~C. Doyle, ``Linear-quadratic-{G}aussian
  heat engines,'' in \emph{Proc. IEEE Conf. Decision Control}, New Orleans, LA,
  2007, pp. 3102--3107.

\bibitem{Hui:ACC:2011}
Q.~Hui, ``Can thermodynamics be used to design control systems?'' in
  \emph{Proc. Amer. Control Conf.}, San Francisco, CA, 2011, pp. 845--850.

\bibitem{HL:Allerton:2011}
Q.~Hui and Z.~Liu, ``A semistabilizability/semidetectability approach to
  semistable $\mathcal{H}_{2}$ and $\mathcal{H}_{\infty}$ control problems,''
  in \emph{49th Ann. Allerton Conf. Comm., Control, Computing}, Monticello, IL,
  2011, pp. 566--571.

\bibitem{BB:ACC:1999}
S.~P. Bhat and D.~S. Bernstein, ``Lyapunov analysis of semistability,'' in
  \emph{Proc. Amer. Control Conf.}, San Diego, CA, 1999, pp. 1608--1612.

\bibitem{HCH:2009}
W.~M. Haddad, V.~Chellaboina, and Q.~Hui, \emph{Nonnegative and Compartmental
  Dynamical Systems}.\hskip 1em plus 0.5em minus 0.4em\relax Princeton, NJ:
  Princeton Univ. Press, 2010.

\bibitem{CHHR:SCL:2008}
V.~Chellaboina, W.~M. Haddad, Q.~Hui, and J.~Ramakrishnan, ``On system state
  equipartitioning and semistability in network dynamical systems with
  arbitrary time-delays,'' \emph{Syst. Control Lett.}, vol.~57, pp. 670--679,
  2008.

\bibitem{HH:AUT:2008}
Q.~Hui and W.~M. Haddad, ``Distributed nonlinear control algorithms for network
  consensus,'' \emph{Automatica}, vol.~44, pp. 2375--2381, 2008.

\bibitem{Hui:TAC:2012}
Q.~Hui, ``Semistability of nonlinear systems having a connected set of
  equilibria and time-delays,'' \emph{IEEE Trans. Autom. Control}, vol.~57,
  no.~10, pp. 2615--2620, 2012.

\bibitem{HCN:2005}
W.~M. Haddad, V.~Chellaboina, and S.~G. Nersesov, \emph{Thermodynamics: A
  Dynamical Systems Approach}.\hskip 1em plus 0.5em minus 0.4em\relax
  Princeton, NJ: Princeton Univ. Press, 2005.

\bibitem{BB:SICON:2003}
S.~P. Bhat and D.~S. Bernstein, ``Nontangency-based {L}yapunov tests for
  convergence and stability in systems having a continuum of equilibra,''
  \emph{SIAM J. Control Optim.}, vol.~42, pp. 1745--1775, 2003.

\bibitem{BB:ACC:2003}
------, ``Arc-length-based {L}yapunov tests for convergence and stability in
  systems having a continuum of equilibria,'' in \emph{Proc. Amer. Control
  Conf.}, Denver, CO, 2003, pp. 2961--2966.

\bibitem{KE:ICNN:1995}
J.~Kennedy and R.~Eberhart, ``Particle swarm optimization,'' in \emph{Proc.
  IEEE Int. Conf. Neural Networks}, Perth, Australia, 1995, pp. 1942--1946.

\bibitem{HL:CDC:2012b}
Q.~Hui and Z.~Liu, ``Semistability-based robust and optimal control design for
  network systems,'' in \emph{Proc. IEEE Conf. Decision Control}, Maui, HI,
  2012, pp. 7049--7054.

\bibitem{HHB:NAHS:2011}
Q.~Hui, W.~M. Haddad, and J.~M. Bailey, ``Multistability, bifurcations, and
  biological neural networks: {A} synaptic drive firing model for cerebral
  cortex transition in the induction of general anesthesia,'' \emph{Nonlinear
  Analysis: Hybrid Systems}, vol.~5, no.~3, pp. 554--572, 2011.

\bibitem{Bernstein:2009}
D.~S. Bernstein, \emph{Matrix Mathematics, \emph{2nd ed.}}\hskip 1em plus 0.5em
  minus 0.4em\relax Princeton, NJ: Princeton Univ. Press, 2009.

\bibitem{ZM:TAC:2011}
D.~Zelazo and M.~Mesbahi, ``Graph-theoretic analysis and synthesis of relative
  sensing networks,'' \emph{IEEE Trans. Autom. Control}, vol.~55, no.~5, pp.
  971--982, 2011.

\bibitem{BP:79}
A.~Berman and R.~J. Plemmons, \emph{Nonnegative {M}atrices in the
  {M}athematical {S}ciences}.\hskip 1em plus 0.5em minus 0.4em\relax New York:
  Academic, 1979.

\bibitem{SZ:SIAP:1970}
J.~Snyders and M.~Zakai, ``On nonnegative solutions of the equation
  {$AD+DA'=C$},'' \emph{SIAM J. Appl. Math.}, vol.~18, pp. 704--714, 1970.

\bibitem{Haddad:TR}
W.~M. Haddad, ``The linear quadratic regulator problem: {A} convex optimization
  approach,'' Sch. Aero. Eng., Georgia Inst. Tech., Atlanta, GA, Tech. Rep.,
  1995.

\bibitem{BH:1995}
D.~S. Bernstein and W.~M. Haddad, \emph{Control-System Synthesis: The
  Fixed-Structure Approach}.\hskip 1em plus 0.5em minus 0.4em\relax Atlanta,
  GA: Monograph, 1995.

\bibitem{probdef}
J.~J. Liang, T.~P. Runarsson, E.~M. Montes, M.~Clerc, P.~N. Suganthan, C.~A.
  {Coello Coello}, and K.~Deb, ``Problem definitions and evaluation criteria
  for the {CEC} 2006 special session on constrained real-parameter
  optimization,'' in \emph{IEEE Conf. Evolutionary Computing}, Vancouver,
  Canada, 2006.

\bibitem{LRMCSCD:TR}
------, ``Problem definitions and evaluation criteria for the {CEC} 2006
  special session on constrained real-parameter optimization,'' Sch. Electric.
  Electron. Eng., Nanyang. Tech. Univ., Singapore, Tech. Rep., 2006.

\bibitem{bestpso}
H.~Lu and W.~Chen, ``Dynamic-objective particle swarm optimization for
  constrained optimizaiton problems,'' \emph{J. Comb. Optim.}, vol.~12, no.~4,
  pp. 409--419, 2006.

\bibitem{bilevel}
X.~Li, P.~Tian, and X.~Min, ``A hierarchical particle swarm optimization for
  solving bilevel programming problems,'' in \emph{Artificial Intelligence and
  Soft Computing -- ICAISC 2006}, ser. Lecture Notes in Computer Science,
  L.~Rutkowski, R.~Tadeusiewicz, L.~A. Zadeh, and J.~M. Zurada, Eds.\hskip 1em
  plus 0.5em minus 0.4em\relax Berlin, Germany: Springer-Verlag, vol.
  4029/2006, pp. 1169--1178.

\bibitem{3rules}
K.~Deb, ``An efficient constraint handling method for genetic algorithms,''
  \emph{Comput. Meth. Appl. Mech. Eng.}, vol. 186, pp. 311--338, 2000.

\bibitem{OFM:PIEEE:2007}
R.~{Olfati-Saber}, J.~A. Fax, and R.~M. Murray, ``Consensus and cooperation in
  networked multi-agent systems,'' \emph{Proc. IEEE}, vol.~95, pp. 215--233,
  2007.

\bibitem{Strogatz:Nature:2001}
S.~H. Strogatz, ``Exploring complex networks,'' \emph{Nature}, vol. 410, pp.
  268--276, 2001.

\bibitem{Hui:ACC:2012}
Q.~Hui, ``A semistability-based design framework for optimal consensus seeking
  of multiagent systems in a noisy environment,'' in \emph{Proc. Amer. Control
  Conf.}, Montr\'{e}al, Canada, 2012, pp. 20--25.

\end{thebibliography}

\end{document}